\documentclass[12pt,reqno]{amsart}

\usepackage[arrow,matrix,curve]{xy}

\usepackage[dvips]{graphicx} 

\usepackage{amssymb, latexsym, amsmath, amscd, array, hyperref,
  epigraph
}

\usepackage{breakurl}   

\usepackage[LGR,T1]{fontenc} 
\usepackage{wrapfig}  
\usepackage{lscape}   
\usepackage{rotating} 

\newcommand{\adscriptiota}{{
  \usefont{LGR}{cmr}{m}{n}\symbol{"7C}%
}}

\DeclareRobustCommand{\ai}[1]{{\mathpalette\doai{#1}}}  
\newcommand{\doai}[2]{%
\oalign{%
~$#1#2$\cr
\hidewidth$#1\text{\,\adscriptiota}$\hidewidth\cr
}%
}

\newcommand\zoon{$\zeta\!\!\ai{\tilde\omega}o\nu$}

\theoremstyle{definition}

\numberwithin{equation}{section}

\numberwithin{equation}{section}

\title[Leibniz on bodies and infinities] {Leibniz on bodies and
  infinities: \emph{rerum natura} and mathematical fictions}

\author[M. Katz]{Mikhail G. Katz} \address{M.~Katz, Department of
  Mathematics, Bar Ilan University, Ramat Gan 5290002 Israel}
\email{katzmik@math.biu.ac.il}

\author[K. Kuhlemann]{Karl Kuhlemann}\address{K. Kuhlemann, Gottfried
  Wilhelm Leibniz University Hannover, D-30167 Hannover, Germany}
\email{kus.kuhlemann@t-online.de}

\author[D. Sherry]{David Sherry} \address{D. Sherry, Department of
  Philosophy, Northern Arizona University, Flagstaff, AZ 86011, US}
\email{David.Sherry@nau.edu}

\author[M. Ugaglia]{Monica Ugaglia} \address{M. Ugaglia, Il Gallo
  Silvestre, Localit\`a Collina 38, Montecassiano, Italy}
\email{monica.ugaglia@gmail.com}

\subjclass[2020]{Primary 01A45       
}

\begin{document}

\thispagestyle{empty}


\keywords{Body; substance; monad; infinitesimal calculus; useful
  fiction; Leibnizian metaphysics; inassignable quantities; infinity;
  infinitesimals; multitude; magnitude; Aristotle; Bernoulli; des
  Bosses; Huygens; Leibniz; Masson; Thomasius; Varignon;
  Wallis}

\begin{abstract}
The way Leibniz applied his philosophy to mathematics has been the
subject of longstanding debates.  A key piece of evidence is his
letter to Masson on bodies.  We offer an interpretation of this often
misunderstood text, dealing with the status of infinite divisibility
in \emph{nature}, rather than in \emph{mathematics}.  In line with
this distinction, we offer a reading of the fictionality of
infinitesimals.  The letter has been claimed to support a reading of
infinitesimals according to which they are logical fictions,
contradictory in their definition, and thus abso\-lutely impossible.
The advocates of such a reading have lumped infinitesimals with
infinite wholes, which are rejected by Leibniz as contradicting the
part-whole principle.  Far from supporting this reading, the letter is
arguably consistent with the view that infinitesimals, as inassignable
quantities, are \emph{mentis fictiones}, i.e., (well-founded) fictions
usable in mathematics, but possibly contrary to the Leibnizian
principle of the harmony of things and not necessarily idealizing
anything in \emph{rerum natura}.  Unlike infinite wholes,
infinitesimals -- as well as imaginary roots and other well-founded
fictions -- may involve accidental (as opposed to absolute)
impossibilities, in accordance with the Leibnizian theories of
knowledge and modality.
\end{abstract}

\maketitle

\today

\tableofcontents

\epigraph{Just as a bounded infinite line is made up of finite ones,
  so a finite line is made up of infinitely small ones.\\ --Leibniz,
  \emph{De Quadratura Arithmetica}}

\epigraph{Infinite and infinitely small quantities can be written out
  of the mathematics altogether via a syncategorematic analysis in
  favour of expressions referring only to finite quantities and their
  relations. --Levey, 2021}

\epigraph{Calculus necessary leads to them, and people who are not
  sufficiently expert in such matters get entangled and think they
  have reached an absurdity.  --Leibniz, \emph{Elementa nova matheseos
    universalis}}

\section
{Toland's indictment of mathematicians}
\label{s1}

%
In 1716, an anonymous critic claimed to identify a confusion shared by
Leibniz and other philosophers of a mathematical bent:
\begin{enumerate}
\item[]
Mais de s'imaginer, qu'ils pourront rendre compte de la \emph{nature
des choses} par de tels Calculs, c'est l\`a precisement que consiste
leur erreur.%
\footnote{Translation: ``But to fancy themselves that they could
account for the nature of things by such Calculations, that is
precisely where their error lies.''}
(\cite[p.\,131]{To11}; emphasis added)
\end{enumerate}
The critic has by now been definitively identified as John Toland%
\footnote{Beeley mistakenly attributes the criticism to Bayle in
  (\cite{Be15} 2015, p.\;26).}
(1670--1722); see Lamarra (\cite{La90}, 1990), Woolhouse (\cite{Wo98},
1998).  Toland accused Leibniz of allowing his calculus to infect his
metaphysics.  The critique was apparently written already in 1703%
\footnote{The date appearing at the end of Toland's text is in 1703,
but there is internal evidence that the piece was touched up in 1714
at the earliest; see \cite{La90}, \cite{Wo98}.}
but only published in a 1716 volume of \emph{Histoire Critique de la
Republique des Lettres}.

In the same year (his last), Leibniz penned a response in a detailed letter%
\footnote{Woolhouse and Francks \cite{Wo97} date the letter 21 august
1716.}
to editor Samuel Masson.  The letter deals with issues of the
philosophy of nature and also comments briefly upon the infinitesimal
calculus.  Bassler (\cite{Ba98}, 1998), Arthur (\cite{Ar19b}, 2019),
and Rabouin and Arthur (RA) (\cite{Ra20}, 2020) have appealed to one
such comment in support of the claim that Leibnizian infinitesimals
are \emph{syncategorematic} (in the sense detailed in
Section~\ref{s4}).  Contextualizing the letter will help evaluate such
claims.

We analyze the Leibnizian exposition on the philosophy of nature and
its historical and theological context in Section~\ref{s2}, where we
also deal with the meaning of his comments on the calculus.  In
Section~\ref{s4} we analyze RA's reading and show that the 1716
comments on the calculus not only provide no support for an
Ishiguro-syncategorematic reading, but support a rather different
interpretation of his infinitesimals: they are well-founded `fictions
of the mind' (\emph{mentis fictiones},%
\footnote{See Section~\ref{s31} for the full quotation.}
Leibniz to des Bosses \cite{Le06}, 1706).  In Section~\ref{s2b} we
analyze the distinctions infinite number \emph{vs} infinite whole, and
bounded infinity \emph{vs} unbounded infinity in Leibniz, as well as
his comparison of the hornangle and inassignables.  In
Section~\ref{s2c} we analyze the meaning of \emph{infinity},
\emph{fiction}, and \emph{well-founded fiction} in Leibniz.  Here we
show that, while infinite wholes contradicting the part-whole
principle are absolute impossibilities, imaginary roots and
infinitesimals are only accidental impossibilities, even if their
definitions are taken to be only nominal.  In Section~\ref{s8} we
present the conclusions of our contextualisation of the letter to
Masson.  The present text extends our analysis of the Leibnizian
heritage pursued in \cite{17b} and elsewhere.

\section
{Monads, \emph{rerum natura}, and mathematics}
\label{s2}

Leibniz opens his letter to Masson with a discussion of the concept
of~{\zoon} (z\=oon),%
\footnote{\label{f6}Gerhardt's edition \cite{Le16} uses the erroneous
  spelling~$\xi\tilde{\omega} o\nu$ (twice) on page 624.  In the
  autograph manuscripts of Leibniz, the word appears
  as~$\zeta\tilde{\omega}o\nu$, without the iota subscript.
%
%
The translation by Ariew \cite[p.\;225]{Le89a} uses the spelling
\emph{zoon}.}
i.e., \emph{living being}.
%
%
%
He describes \emph{monad} as the underlying \emph{substance}%
\footnote{The term \emph{substance} is used here in the technical
sense of a metaphysical substratum of physical beings.  See further in
note~\ref{f8}.}
of the z\=oon.  According to Leibniz, the relation between nature and
monadic reality is such that nature is a phenomenon, whereas the
\emph{monads} are the substances underlying natural phenomena.
Leibniz first used the term `monad' in 1696; in earlier texts he used
various expressions such as `individual substance'
\cite[p.\;32]{Le89a}, `atom of substance', and `metaphysical point'
\cite[p.\,142]{Le89a}.%
\footnote{\label{f8}Early work by Leibniz on substance dating from
  1668 aimed ``to effect a reconciliation between Roman Catholics and
  Protestants.  {\ldots} These works are especially valuable for what
  they reveal about the motivations behind Leibniz's first account of
  substance'' (Mercer--Sleigh \cite[p.\;68]{Me94}).  In a letter to
  des Bosses dated 8\;september 1709, Leibniz distances himself from
  both transubstantiation and consubstantiation, and sketches a
  monad-based approach \cite[p.\,153]{Le07}.}

\subsection{Leibnizian Passage I}
\label{s21}

Toward the end of the letter, Leibniz turns to Toland's claim that
mathematicians are not successful as philosophers.  After a bit of ad
hominem, Leibniz gets down to the business of refuting the claim.  His
strategy is to draw a line between
\begin{enumerate}
\item
philosophy, concerned with the `nature des choses', i.e., \emph{rerum
  natura}; and
\item
mathematics, concerned with applying entities, both ideal and
fictional, in geometry and physics.
\end{enumerate}
Toland accused Leibniz of viewing the extension (i.e., the continuum)
as made up of mathematical points (the punctiform view).%
\footnote{Toland's criticism of Leibniz being so widely off the mark
  has led Stuart Brown to speculate that Toland ``mistakenly supposed
  that Leibniz was cast in the same philosophical mould as [Joseph]
  Raphson'' \cite[note\;51]{Br99}, and copied over some criticisms
  targeting Raphson from another text of Toland's.  Raphson was an
  associate of Newton's.}
Leibniz first responds to this accusation, and then comments on the
relation between the calculus and the \emph{nature des choses}:
\begin{enumerate}
\item[] [\textbf{Passage I}:] I am also far from \emph{making
  extension up of mathematical points}.~\,{\ldots} And,
  notwithstanding my \emph{infinitesimal calculus}, I do not admit any
  real [\emph{veritable}] infinite number, even though I confess that
  the multitude of things surpasses any finite number [\emph{la
      multitude des choses passe tout nombre fini}], or rather any
  number.  (Leibniz \cite{Le16} as translated by Ariew in
  \cite[p.\;229]{Le89a})
\end{enumerate}
Passage\;I draws a line between, on the one hand, the mathematician's
task of exploiting the well-founded fictions of the infinitesimal
calculus, and on the other, the philosopher's task of elucidating the
natural phenomena (and perhaps the ultimately real monadic entities
which underlie the phenomena) in the framework of an unequivocal
rejection of infinite wholes.

In Sections~\ref{anto} and \ref{anto2} we will deal with the
distinction Mathematics \emph{vs} \emph{rerum natura}.  The
distinction infinite wholes \emph{vs} well-founded fictions will be
analyzed in Section~\ref{s2b}.

\subsection{Leibnizian Passage II}
\label{anto}

Unlike the \emph{choses} whose nature Leibniz seeks to explore in the
1716 letter, mathematical entities exploited in the infinitesimal
calculus are only useful fictions (\emph{mentis fictiones}; see
Section~\ref{s31}):
\begin{enumerate}\item[]
[\textbf{Passage II}:] The infinitesimal calculus is useful with
respect to the application of mathematics to physics; however, that is
not how I claim to account for the nature of things [\emph{la nature
    des choses}].  For I consider infinitesimal quantities to be
useful fictions.  (Leibniz \cite{Le16} as translated by Ariew in
\cite{Le89a}, p.\;230)
\end{enumerate}
Leibniz's comment ``however, that is not how I claim to account for
the nature of things'' is a direct response to Toland's allegation
that Leibniz seeks to account for the \emph{nature des choses} by
means of his calculus (see Section~\ref{s1}).  Garber quotes
Passage~II and notes:
\begin{enumerate}\item[]
[T]he point seems to be that nature is one thing, and its mathematical
representation is another.  \cite[p.\;303]{Ga08}
\end{enumerate}
Garber points out a significant difference between the positions of
Leibniz and the Cartesians:
\begin{enumerate}\item[]
[Leibniz's] opponents are the Cartesians who have tried to make nature
mathematical in a literal sense, to make the physical world over into
a physical instantiation of mathematical concepts [whereas Leibniz]
can embrace the mathematical representation of dead force in terms of
infinitesimals, without having to say that there are real
infinitesimals in nature.  \cite[p.\;306]{Ga08}
\end{enumerate}
Thus in Garber's view, in Passage II Leibniz (disagreeing with the
Cartesians) insists on the separation of \emph{rerum natura} and its
mathematical representation.

\subsection
{Leibnizian evolution on mathematics and \emph{rerum natura}}
\label{anto2}

The distinction between the mathematical realm and the \emph{rerum
  natura} is a crucial feature of Leibniz's mature philosophy.
Whereas he started with a belief that physics could be reduced to
mechanics,%
\footnote{See the 1668--9 correspondence with Thomasius, Leibniz
  \cite{Le23} A2.1$^2$.\;16--44; cf.\;\cite[p.\;71]{Me94}.}
and hence to mathematics, over the years mathematics ceased to be
perceived by Leibniz as the foundation of physics and turned into a
mere representation thereof.

For the young Leibniz, extension (\emph{corpus mathematicum}) is a
component of the matter of \emph{choses}.  As Leibniz explains to
Thomasius, a physical body is the compound of matter -- which Leibniz
identifies with extension and impenetrability -- and form, which he
identifies with shape:
\begin{enumerate}\item[]
Space is a primary extended being or a mathematical body (\emph{corpus
  mathematicum}), which contains nothing but three dimensions and is
the universal locus of all things.  Matter is a secondary extended
being, or that which has, in addition to extension or mathematical
body, also a physical body (\emph{corpus physicum}), that is,
resistance, antitypy, solidity, the property of filling space,
impenetrability.%
\footnote{Resistance, antitypy, solidity, the property of filling
  space and impenetrability are synonyms, and refer to the
  impossibility, for a physical body, of being in the same space with
  another thing, as Leibniz explains a few lines later.}
(Leibniz to Thomasius, 20--30 april 1669,
\cite[p.\;34]{Le69}. Transl.\;Loemker \cite[p.\,100]{Le89} with minor
changes)%
\footnote{\;``Spatium est Ens primo-extensum, seu corpus mathematicum,
  quod scilicet nihil aliud continet quam tres dimensiones, estque
  locus ille universalis omnium rerum. Materia est ens
  secundo-extensum, seu quod praeter extensionem vel corpus
  mathematicum habet et corpus physicum, id est, resistentiam,\,
  '\!\!\!\!$\alpha\nu\tau\iota\tau\upsilon\pi\acute\iota\alpha\nu$,
%
%
  crassitiem, spatii-repletivitatem, impenetrabilitatem.''}
\end{enumerate}
During this phase of Leibniz's development, physical bodies are seen
as continuous, exactly as extension is.  Namely, they are potentially
infinitely divisible:
\begin{enumerate}\item[]
Matter has quantity too, though it is indefinite, or interminate as
the Averroists call it.  For being continuous, it is not cut into
parts and therefore does not actually have boundaries.  (Leibniz to
Thomasius, 20--30 april 1669,
\cite[p.\;26--27]{Le69}. Transl.\;Loemker \cite[p.\;95]{Le89})%
\footnote{\;``Quantitatem quoque habet materia, sed interminatam, ut
  vocant Averroistae, seu indefinitam, dum enim continua est, in
  partes secta non est, ergo nec termini in ea actu dantur.''}
\end{enumerate} 
By contrast, in his mature system, Leibniz distinguishes between
mathematical extension, which is \emph{potentially} infinitely
divisible, and the matter of physical bodies, which is \emph{actually}
infinitely divided:
\begin{enumerate}\item[]
But in real things, that is, bodies, the parts are not indefinite --
as they are in space, which is a mental thing -- but actually
specified in a fixed way according to the divisions and subdivisions
which nature actually introduces through the varieties of motion. And
granted that these divisions proceed to infinity, they are
nevertheless the result of fixed primary constituents or real unities,
though infinite in number.  Accurately speaking, however, matter is
not composed of these constitutive unities but results from them.
(Leibniz to de Volder, 30\;june 1704 \cite[p.\;268]{Le04}.
Transl.\;Loemker \cite[p.\;536]{Le89})%
\footnote{\;``At in realibus, nempe corporibus, partes non sunt
  indefinitae (ut in spatio, re mentali), sed actu assignatae certo
  modo, prout natura divisiones et subdivisiones actu secundum motuum
  varietates instituit, et licet eae divisiones procedant in
  infinitum, non ideo tamen minus omnia resultant ex certis primis
  constitutivis seu unitatibus realibus, sed numero
  infinitis. Accurate autem loquendo materia non componitur ex
  unitatibus constitutivis, sed ex iis resultat.''}
\end{enumerate}
Similarly in a 11 march 1706 letter to des Bosses, Leibniz writes:
\begin{enumerate}\item[]
To pass now from the ideas of geometry to the realities of physics, I
hold that matter is actually fragmented into parts smaller than any
given, or that there is no part of matter that is not actually
subdivided into others exercising different motions. This is demanded
by the nature of matter and motion and by the structure of the
universe, for physical, mathematical, and metaphysical reasons.
(Leibniz to des Bosses, 11\;march 1706 \cite[p.\;305]{Le06}.
Transl.\;Look--Rutherford \cite[p.\;33--35]{Le07})%
\footnote{\;``Caeterum ut ab ideis Geometriae, ad realia Physicae
  transeam; statuo materiam actu fractam esse in partes quavis data
  minores, seu nullam esse partem, quae non actu in alias sit
  subdivisa diversos motus exercentes. Id postulat natura materiae et
  motus, et tota rerum compages, per physicas, mathematicas et
  metaphysicas rationes.''}
\end{enumerate}
The last passage is analyzed by Antognazza, who emphasizes the
distinction in Leibnizian thought between the mathematical realm and
\emph{rerum natura}:
\begin{enumerate}\item[]
Leibniz is quite consistent in pointing out that the actual infinite
he is endorsing concerns the `real' as opposed to the `ideal' order.
In the letter of 11 March 1706 to des Bosses {\ldots} he explicitly
stresses that in moving his attention to the actual infinite, he is
shifting from the ideal to the real order.
%
%
(Antognazza \cite{An15}, 2015, p.\;9)
\end{enumerate}
To avoid confusion between continuous extension, which is an ideal
entity, and the actual structure of physical matter, Leibniz employs
for the latter the adjective \emph{contiguous}:
\begin{enumerate}\item[]
I recall that Aristotle too distinguishes between Contiguum and
Continuum: things are continuous if their extremes are one, and are
contiguous if their extremes are together.  (\emph{Pacidius
  Philaleti}, 1676 (Leibniz \cite{Le23}.  A6.3. 537)%
\footnote{\;``Memini Aristotelem quoque Contiguum a Continuo ita
  discernere, ut Continua sint quorum extrema unum sunt, Contingua
  quorum extrema simul sunt.''  Cf. Aristotle Physics VI.1 231a
  22--3.}
\end{enumerate}
Levey similarly emphasizes the difference between the actual
subdivisions of the real world and the potential ones of the ideal
mathematical world:
\begin{enumerate}\item[]
[To Leibniz,] a body is separable into various parts because it
\emph{actually} has contiguous parts that cohere together but which
could be brought not to cohere and be separated from one
another. {\ldots} Potentiality, in the sense of potential divisions or
potential parts, is a concept that belongs to the `ideal' realm of
mathematics and geometry but has no application to the world of
matter.  (Levey \cite{Le98}, 1998, p.\;53, note\;6; emphasis in the
original)
\end{enumerate}
Views similar to that of Antognazza were expressed by Bosinelli
\cite[p.\,168]{Bo91} and Breger \cite[p.\,124]{Br16}, as acknowledged
by Arthur in \cite[p.\,156]{Ar18}.
%
%
But especially, they were expressed by Leibniz himself, who
specifically warned his readers against the misunderstandings arising
from the conflation of the two realms:
\begin{enumerate}\item[]
It is the confusion of the ideal with the actual which has muddled
everything and caused the labyrinth of the composition of the
continuum.  (Remarques sur les Objections de M. Foucher, 1695.
Gerhardt \cite{Ge49} vol IV, p.\;491. Transl.\;Ariew--Garber
\cite[p.\,146]{Le89a})%
\footnote{\;``Et c'est la confusion de l'ideal et de l'actuel qui a tout
  embrouill\'e et fait le labyrinthe de \emph{compositione
    continui}.''}
\end{enumerate}
A detailed study of the issue appears in Ugaglia (\cite{Ug22}, 2022).

Just as the Leibnizian Passage\;I in Section~\ref{s21}, Passage\;II in
Section~\ref{anto} draws a line between the fictional entities
(\emph{mentis fictiones}) of the infinitesimal calculus --
well-founded fictions, or useful fictions, as Leibniz also calls them
-- and the entities of \emph{rerum natura}.  The conflation of the two
realms is behind some of the purported evidence in favor of the
so-called syncategorematic reading, analyzed in Section~\ref{s4}.

\section
{Syncategorematics \emph{vs} keeping bodies and calculus separate}
\label{s4}

RA follow Ishiguro (\cite{Is90}, 1990, Chapter\;5) in interpreting the
Leibnizian term `useful fiction' \emph{syncategorematically} via a
logical analysis of the quantifier clause `smaller than every given'.
This stance enables RA to claim that Leibnizian infinitesimals are `in
keeping with the Archime\-dean axiom' \cite[Abstract]{Ra20} (in the
previous installment by Arthur \cite{Ar13} the refrain was `fully in
accord with the Archimedean Axiom').

\subsection
{Small magnitudes \emph{vs} Ishiguro's alternating quantifiers}
\label{s41}

In 1990, Ishiguro formulated the following hypothesis concerning
Leibnizian infinitesimals:
\begin{enumerate}\item[]
It seems then that throughout his working life as a mathematician
Leibniz did not think of founding the calculus in terms of a
\emph{special kind of small magnitude}.  \cite[p.\;86]{Is90} (emphasis
added)
\end{enumerate}
To explain in what sense Leibniz allegedly did not found the calculus
on a special kind of small magnitude, she elaborates as follows:
\begin{enumerate}\item[]
  It seems that when we make reference to infinitesimals in a
  proposition, we are not designating a fixed magnitude incomparably
  smaller than our \emph{ordinary} magnitudes.  Leibniz is saying that
  whatever small magnitude an opponent may present, one can assert the
  existence of a smaller magnitude.  In other words, we can paraphrase
  the proposition with a \emph{universal} proposition with an embedded
  \emph{existential} claim.  (ibid., p.\;87; emphasis on `ordinary',
  `universal', and `existential' added)
\end{enumerate}
Accordingly, when Leibniz asserted that his inassignable~$dx$, or
alternatively~$\epsilon$, was smaller than every assignable
quantity~$Q$, he really meant that for each given~$Q>0$ there exists
an `ordinary' \mbox{$\epsilon>0$} such that~\mbox{$\epsilon<Q$}.

Ishiguro interprets Leibniz's uses of the term `infinite quantity' by
a similar paraphrase involving embedded quantifiers.  Such an approach
has been described as a \emph{syncategorematic} treatment of
statements involving `infinite' and its cognates.  Ishiguro goes on to
compare the theories of Leibniz and Russell and concludes:
\begin{enumerate}\item[]
The similarity which we find in the two theories lies in the common
intention Leibniz and Russell have of understanding a sentence which
ostensibly designates a specific entity, as really, in its
\emph{logical} form, being a quantified sentence; i.e., a universal or
an existential sentence.  (ibid., p.\;99; emphasis on `logical' added)
\end{enumerate}
Taking her cue from Russell, Ishiguro reads Leibnizian `useful
fictions' as \emph{logical fictions}.%
\footnote{An alternative reading of Leibnizian useful fictions as
  \emph{pure fictions} was developed in \cite{14c}.  For more details
  see \cite{16a}.}
In the same vein, Levey writes:
\begin{enumerate}\item[]
[T]he fiction is intended by Leibniz to be `logical' in character:
infinite and infinitely small quantities can be written out of the
mathematics altogether via a syncategorematic analysis in favour of
expressions referring only to finite quantities and their relations.
\cite[p.\,148]{Le21}
\end{enumerate}
By a similar quantifier sleight-of-hand, RA declare Leibniz's uses of
the term `infinitesimal' to be ``fully-in-accord/in-keeping with the
Archi\-me\-dean axiom'' (see further in Section~\ref{s43}).

\subsection
{Rabouin on syncategorematic entities}
\label{s42}

In 2013, Rabouin toys with ``the idea that the entities being studied
are \emph{relational} or in Leibniz's parlance `syncategorematic'{}''
\cite[p.\,120]{Ra13} and adds: ``This is the reason why Leibniz calls
them \emph{fictitia} (they are terms not referring to individual
beings, but to some relational properties)'' (ad loc., note 40).

On the other hand, in 2015, Rabouin seeks to distance himself from
attributing hidden quantifiers to Leibniz, in the following terms:
\begin{enumerate}\item[]
The core of the demonstration {\ldots} is the \emph{arbitrariness} of
the choice of~$\varepsilon$.  But this arbitrariness does not amount,
in modern terms, to a universal quantification (at least in classical
first order logic), which would be meaningless to Leibniz.  (Rabouin
\cite{Ra15}, 2015, p.\;362)
\end{enumerate}
Yet in note\;25 on the same page Rabouin appears to endorse Ishiguro's
reading: ``It is, in Leibniz' terminology, a `syncategorematic' entity%
\footnote{No punctuation mark occurs between `entity' and `see' in
  Rabouin \cite[note\;25]{Ra15}.}
see Ishiguro (1990) {\ldots}'' (ibid.).  Rabouin goes on to claim the
following:
\begin{enumerate}\item[]
  It should then be clear why infinitesimals were called by Leibniz
  `fictions'. In and of itself, there is no such thing as a `quantity
  smaller than any other quantity'.  This would amount to the
  existence of a minimal quantity and one can show that a minimal
  quantity implies \emph{contradiction}. (ibid.; emphasis added)
\end{enumerate}
What are we to make of Rabouin's claimed detection of contradiction?
While it is true that there is no such thing as a nonzero quantity
smaller than any other quantity, Leibniz requires his inassignable
infinitesimals to be smaller, not \emph{than every quantity}, but
rather \emph{than every assignable quantity}; for an illustration in
terms of hornangles see Section~\ref{s53}.  This accords with a
glossary entry in (Rabouin's 2020 coauthor) Arthur's volume of
translations.  Here an infinitesimal is defined as
\begin{enumerate}\item[]
a part smaller than any assignable (see {\tiny INASSIGNABILIS}), a
definition to which Leibniz frequently has recourse.%
\footnote{\label{f19}This entry from Arthur's glossary indicates that
he is aware of the Leibnizian assignable/inassignable dichotomy; see
main text at note~\ref{f30}.}
(Arthur \cite{Ar01}, 2001, p.\;452)
\end{enumerate}
Thus a particular infinitesimal~$\epsilon>0$ will satisfy
\begin{equation}
\label{e51}
\epsilon<\frac12,\; \epsilon<\frac13,\; \epsilon<\frac14,\;
\text{etc.},
\end{equation}
so long as the denominator is assignable.  Due to what appears to be a
mathematical misunderstanding, Rabouin is led to a conclusion that a
bona fide infinitesimal would be contradictory and hence `fictional'.

Once one realizes that the contradiction is not there to begin with,
there is no compelling reason to interpret Leibnizian fictionalism as
the counterpart of an allegedly contradictory nature of
infinitesimals, as Rabouin does in 2015.  A geometric illustration of
\eqref{e51} in terms of hornangles appears in Section~\ref{s53}.

\subsection{Letter to Masson and RA's quantifiers}  
\label{s43}

In Section~\ref{s42}, we examined Rabouin's 2015 attempt to declare
Leibnizian infinitesimals contradictory.  Five years later in 2020,
one finds a related attempt to declare Leibnizian infinitesimals
contradictory by Rabouin and Arthur (RA) in \cite{Ra20} (see also
Section~\ref{s31z}).  Infinitesimals do contradict the Archimedean
axiom, but if this is what RA mean by their contradictory claim, then
their argument (in support of the contention that Leibnizian allegedly
syncategorematic infinitesimals are `in keeping with the Archime\-dean
axiom') would be circular: the positing of the Archime\-dean axiom
predetermines the conclusion that non-Archimedean infinitesimals are
contradictory.%
\footnote{Cantor's published proof that infinitesimals are
  contradictory suffers from a related circularity; see Ehrlich
  (\cite{Eh06}, 2006).
%
%
}

In search of evidence in support of an Ishiguroan interpretation, RA
turn to the 1716 letter to Masson where Leibniz employs the quantifier
clause
\begin{enumerate}\item[]
`la multitude des choses passe tout nombre fini'.
\end{enumerate}
Via the quantifier clause, RA seek to establish a connection between
useful fictions and
%
%
logical fictions, `in keeping with the Archimedean axiom'.

However, in his 1716 letter Leibniz is not referring to mathematical
entities, whether fictional or ideal, when he employs the quantifier
clause above.  Rather, he is referring to natural phenomena and the
monads which underlie them (see Section~\ref{s2}).  Neither natural
phenomena nor the monads are \emph{useful fictions}, unlike
infinitesimals.  The quantifier clause in the 1716 letter therefore
does not refer to mathematical entities such as infinitesimals and
infinite quantities.  Rather, it merely reasserts Leibniz's opposition
to infinite wholes as contradictory.  RA claim the following
similarity:

\begin{enumerate}\item[]
Leibniz insisted that since every body is actually divided by motions
within it into further bodies that are themselves similarly divided
without bound, bodies `are actually infinite,
%
%
that is to say, more bodies can be found than there are unities in any
given number whatever' (A\;VI 4, 1393; LLC 235).  That is, their
multiplicity `surpasses every finite number'.$^{26}$

\quad\emph{Similarly}, when Pierre Varignon asked Leibniz for
clarification of his views on the infinitely small in 1702, Leibniz
replied (February 2, 1702) that it had not been his intention `to
assert that there are in nature infinitely small lines in all rigour,
or compared with ours, nor that there are lines infinitely greater
than ours'.  (\cite[p.\;412]{Ra20}; emphasis on `Similarly' added)
\end{enumerate}
Here RA's footnote 26 contains the following text: 
\begin{enumerate}\item[]
``26 Cf.\;what Leibniz wrote to Samuel Masson in the last year of his
  life: `Notwithstanding my Infinitesimal Calculus, I do not at all
  admit a genuine infinite number, although I confess that the
  multiplicity of things surpasses every finite number, or rather,
  every number.' (GP VI 629)''
\end{enumerate}
We can agree with RA's comments on Leibniz's understanding of actual
subdivision of bodies.  However, RA also claim a \emph{similarity} in
reference to the following Leibnizian texts:
\begin{enumerate}
\item
his \emph{Actu infinitae sunt creaturae} (1678--81);
\item
letter to Varignon (1702);
\item
letter to Masson (1716).
\end{enumerate}
The alleged similarity is apparently based on Leibniz's calling into
question certain notions of infinity in these texts.  We analyze their
claim in Section~\ref{s54}.

\subsection{Comparison of letters to Varignon and Masson}
\label{s54}

In connection with the comparison of the letters, we note the
following three points.

\medskip
\textbf{1.} Notwithstanding the fact that Leibniz briefly mentions the
infinitesimal calculus in the 1716 passage, the substance of the
sentence is not concerned with infinitesimal calculus, but rather is a
reflection on the philosophy of nature that Leibniz deliberately
\emph{contrasted} with the calculus.  What is put into question in the
1716 letter is the concept of an infinite taken as a whole (deemed
contradictory by Leibniz), in the context of a metaphysical analysis
of body and substance.

Meanwhile, the 1702 comment addressed to Varignon deals, not with
metaphysics, but with mathematical fictions (the relevant passage from
the letter to Varignon appears in Section~\ref{s33} at
note~\ref{f49}).  Notice that, unlike the 1716 passage, the 1702
passage does not even mention `number', but speaks rather of geometric
objects, such as infinitely large \emph{lines} (i.e, line segments),
which are instances of Leibnizian \emph{bounded infinity} contrasted
with unbounded infinite wholes (see Sections~\ref{s32}
and~\ref{s34b}).  The 1702 passage reasserts the Leibnizian position
that infinitesimal lines need not idealize anything in nature to be
useful, and in this sense can be taken to be fictional; see
Section~\ref{s2b}.  Therefore RA's claimed similarity has no basis.

\medskip
\textbf{2.} A similar conflation affects RA's use of the Leibnizian
phrase
\begin{enumerate}\item[]
``So bodies are actually infinite, i.e., more bodies can be found than
there are unities in any given number'' (Leibniz as translated by
Arthur in \cite[p.\;235]{Ar01})
\end{enumerate}
part of which is quoted in the passage from RA appearing in
Section~\ref{s43}.  The phrase concerns \emph{bodies} (in \emph{rerum
  natura}) and not infinitesimals.  The distinction was emphasized by
both Antognazza and Levey (see Section~\ref{anto}).

\medskip
\textbf{3.} Referring to possible connections between mathematics and
philosophy in Leibniz and the interpretation that puts differential
calculus at the core of Leibniz' philosophy, Rabouin solo writes:
\begin{enumerate}\item[]
Leibniz was \dots\,very explicit about some connections which he
\emph{resisted making}---although modern commentators tend to put a
lot of emphasis on them.  (\cite{Ra15b}, 2015, note 47, p.\;69;
emphasis on `resisted making' ours)
\end{enumerate}
Rabouin goes on to quote Passage II from the 1716 letter (see
Section~\ref{anto}) as evidence (for a reluctance on the part of
Leibniz to make certain connections between mathematics and
philosophy).  Rabouin's interpretation of Passage\;II, insisting on a
segregation of mathematics and philosophy, is compatible with our
view, but not with the approach pursued by Rabouin and Arthur in
\cite[p.\;412]{Ra20} that seeks to enlist the 1716 letter as evidence
in favor of a syncategorematic reading of the Leibnizian calculus.
%
%
In Sections~\ref{s2b} and \ref{s2c}, we develop an interpretation of
Leibnizian infinitesimals that is more faithful to the Leibnizian
texts.

\section{Well-founded fictions: the mathematics}
\label{s2b}

Leibniz described the infinitesimals of his calculus as
\emph{well-founded} because they perform successfully within the
system of rules he developed for the calculus.  Similarly, imaginaries
are well-founded because mathematical experience shows that they are
useful in the solution of the cubic equation and other problems.

Leibniz's choice of the term \emph{multitude} in Passage I
(Section~\ref{s21}) is significant.  Here Leibniz reinforces his
rejection of a punctiform conception of extension by reminding the
reader that he rejects a real infinite number, meaning that he rejects
an infinite multitude taken as a whole (implied by the punctiform view
of the continuum) as contradictory.

\subsection{Part-whole principle}
\label{s31b}

Leibniz held that infinite wholes would contradict the part-whole
principle.  Already in 1672 Leibniz wrote:
\begin{enumerate}
\item[]
[T]here are as many squares as numbers, that is to say, there are as
many square numbers as there are numbers in the universe.  Which is
impossible.  Hence it follows either that \emph{in the infinite the
whole is not greater than the part}, which is the opinion of Galileo
and Gregory of St.\;Vincent, and \emph{which I cannot accept}; or that
infinity itself is nothing, i.e.\ that it is not one and not a whole.
Or perhaps we should say, distinguishing among infinities, that the
most infinite, i.e.\ all the numbers is something that implies a
contradiction, {\dots} (Leibniz as translated by Arthur in
\cite[p.\;9]{Ar01}; emphasis added)
\end{enumerate}
Thus already in 1672, Leibniz held that the infinity of all numbers --
that we may today refer to as an infinite cardinality%
\footnote{\label{f21}As noted correctly by Ishiguro, ``Leibniz did not
  think that there should be what we call the cardinality of the set
  of all things'' \cite[p.\;80]{Is90}.}
-- contradicts the part-whole principle.%
\footnote{Commentators who wish to lump infinite wholes with
  infinitesimals encounter an immediate difficulty: if their fates
  were bound together in Leibniz's mind, why did he reject infinite
  wholes as contradictory already in 1672 but waited until 1676 (see
  note~\ref{f32} for Arthur's tight timeframe) to rule on
  infinitesimals?}
To explain Leibniz' rejection of infinite wholes in modern terms, one
could perhaps surmise that Leibniz would have rejected as incoherent
the modern notion of infinite cardinality or cardinal number (a point
painstakingly argued in \cite{Ar19b}).
%
%

While Leibniz himself of course does not distinguish between infinite
number and infinite cardinality, he does mention a related distinction
between infinite magnitude and infinite multitude in a letter to
des\;Bosses \cite[p.\;31]{Le06}.  The correspondence between Leibniz
and Bernoulli indicates that Leibniz was clearly aware of the
difference between infinite multitude and infinite number.%
\footnote{\label{f10}Thus, in a rebuttal of Bernoulli's argument for
  the existence of infinitesimals based on an analysis of a geometric
  series, Leibniz points out: ``I will concede the existence of the
  \emph{infinite multitude}, but this multitude is neither a
  \emph{number} nor a coherent \emph{whole}.  It means nothing more
  than that there are more parts, which could be referred to by any
  number at all, just as there is a multitude {\ldots}\;of all
  numbers; this multitude, however, is itself neither a number nor a
  coherent whole'' (Leibniz as translated by Sierksma in \cite{Si99},
  1999, p.\;447).}
The distinction is important because Leibniz’s rejection of an
infinite multitude, or collection, does not entail a rejection of
infinite magnitude, or quantity, or number, provided they are not
considered as a whole.  What leads to contradiction is not infinity in
itself but an infinity taken as a whole.  Regrettably, the notions of
infinite whole and infinite quantity have been conflated in recent
literature.  In particular, such a conflation is at the root of the
Ishiguro-syncategorematic interpretation of infinitesimals, analyzed
in Section~\ref{s4}.

\subsection{Infinite wholes and infinite numbers}
\label{s31z}

Both Arthur (\cite{Ar19c}, 2019) and Rabouin and Arthur (\cite{Ra20},
2020) seek to connect the Leibnizian rejection of infinite wholes and
his description of infinitesimals as \emph{fictional}, and to document
Leibnizian rejection of infinite number.  Thereby they seek to
assimilate infinitesimals to infinite wholes.%
\footnote{\label{f24b}Leibniz rejected both \emph{minima} and
  \emph{maxima} for continua in the following terms:
  ``\emph{Scholium}.  We therefore hold that two things are excluded
  from the realm of intelligibles: minimum, and maximum; the
  indivisible, or what is entirely \emph{one}, and \emph{everything};
  what lacks parts, and what cannot be part of another'' (Leibniz as
  translated by Arthur in \cite[p.\,13]{Ar01}).  The rejection of
  \emph{maxima} is the rejection of infinite wholes (e.g., unbounded
  lines).  The rejection of their counterparts, \emph{minima}, is the
  rejection of putative simplest constituents of the continuum, i.e.,
  the rejection of a punctiform continuum (see Section~\ref{s21}).  To
  Leibniz, points play only the role of endpoints of line segments.
  Thus the rejected counterparts of the contradictory infinite wholes
  are, not infinitesimals as per RA, but rather points viewed as the
  simplest constituents of a continuum.}
RA claim the following:
\begin{enumerate}
\item[] Leibniz never changed his mind and always claimed that
  \emph{such entities} do not exist because they imply a
  contradiction.  One can nonetheless use infinities and infinitely
  small quantities in a calculation provided one can furnish a way of
  doing correct demonstrations with them; that is, as long as one can
  identify \emph{conditions} under which their use will not lead to
  error.  (\cite[p.\;413]{Ra20}; emphasis on `such entities' and
  `conditions' added)
\end{enumerate}
While no such conditions are ever identified by RA, in their
Conclusion they speak of using infinitesimals `under certain
specified conditions' \cite[p.\;441]{Ra20}, again unspecified.%
\footnote{\label{f24}See also note~\ref{f48} on well-foundedness and
contradictions.}
The key issue, however, is what is meant exactly by ``such entities''
mentioned by RA.  While Leibniz rejected as contradictory an infinite
multitude taken as a whole (see Section~\ref{s2}), RA provide no
evidence that Leibniz viewed infinitesimal and infinite quantities as
contradictory (and not merely as well-founded fictions).%
\footnote{\label{f25}RA conflate multitudes and quantities when they
claim that
``[Leibniz] held that the part-whole axiom is constitutive of
quantity, so that the concept of an infinite quantity, such as an
infinite number or an infinite whole, involves a contradiction''
\cite[p.\;406]{Ra20}.  The part-whole axiom is constitutive of
multitude but not of quantity.  As we mentioned in note~\ref{f10},
Leibniz clearly understood the difference between infinite multitude
and infinite quantity.}

A similar conflation appears in Rabouin solo, who claims that 
\begin{enumerate}\item[]
[Leibniz] regularly stumbled upon the fact that an `infinite number'
or multitude was a notion entailing a contradiction.  Infinitely small
numbers, of the kind introduced by Wallis in his \emph{Arithmetics of
  the infinites} \emph{(Arithmetica infinitorum}) as inverse of
infinite number, naturally inherited this qualification.  (Rabouin
\cite{Ra21}, 2021)
\end{enumerate}
The contradictory nature of infinite wholes is `inherited' not by
infinitesimals but by \emph{minima}, i.e., points viewed as the
simplest constituents of the continuum.%
\footnote{RA's confusion of \emph{minima} and infinitesimals was
  already noted in note~\ref{f24b}.}

Arthur similarly conflates infinite collections and infinite number
when he claims the following:
\begin{enumerate}
\item[]
[Leibniz] held that there are actually infinitely many substances,
actually infinitely many parts of matter, and actually infinitely many
terms in an infinite series:%
\footnote{Arthur has never given any evidence for his claim that
  Leibniz viewed the terms of an infinite series as an actual infinity
  (even distributively).  Antognazza notes that Leibniz ``is offering
  a mathematical analogy, as opposed to maintaining that the actual
  infinite (even if thought of syncategorematically) applies to
  mathematical, abstract entities, and to the ideal, mathematical
  continuum'' \cite[p.\;9]{An15}.  See Section~\ref{anto2} for further
  details.}
that is, there are so many that, however many are assigned, there are
more, but there is no infinite collection of them, and
\emph{therefore} no infinite number.  (Arthur \cite{Ar19c}, 2019,
p.\,152; emphasis added)
\end{enumerate}
Arthur's final `therefore' is a non-sequitur if the expression
`infinite number' is meant to include infinite quantity (in addition
to infinite wholes).


The related issue of nominal definitions is discussed in Section
\ref{s47}.

\subsection{Bounded and unbounded infinity}
\label{s32}

A Leibnizian distinction of long standing -- at least since his
\emph{De Quadra\-tura Arithmetica} (DQA), Propositio\;XI and the
Scholium following it -- is between bounded infinity and unbounded
infinity; see e.g., Knobloch \cite[p.\;42]{Kn90},
\cite[pp.\;266--267]{Kn94}.
%
%
In this text of fundamental importance for the foundations of the
calculus, Leibniz contrasts bounded infinity and unbounded infinity in
the following terms:
\begin{enumerate}\item[]
But as far as the activity of the mind%
\footnote{For Leibniz on minds, see Section~\ref{s31} and main text
at note~\ref{f35}.}
with which we measure infinite areas is concerned, it contains nothing
unusual because it is based on a certain fiction and proceeds
effortlessly on the assumption of a certain, though \emph{bounded, but
infinite line}; therefore it has no greater difficulty than if we were
to measure an area that is finite in length.  {\ldots} Just as points,
even of infinite numbers, are unsuccessfully added to and subtracted
from a bounded line, so a bounded line can neither form nor exhaust an
unbounded one, however many times it has been repeated.  This is
different with a \emph{bounded but infinite line} thought to be
created by some multitude of finite lines, although this multitude
exceeds any number.  And just as a \emph{bounded infinite line} is
made up of finite ones, so a finite line is made up of infinitely
small ones, yet divisible.  (Leibniz, DQA \cite{Le04b}, Scholium
following Propositio\;XI; translation ours; emphasis added)
\end{enumerate}
It is worth noting that bounded/unbounded infinity is not the same
distinction as potential/actual infinity.  While actual infinity
(understood distributively) is possible in the material realm (see
Section~\ref{anto2}), unbounded infinity (understood collectively,
i.e., as a whole) is a contradictory concept to Leibniz.  Bounded
infinity is a term Leibniz reserves mainly to discuss the well-founded
fictions used in his infinitesimal calculus, namely the infinitely
large and (its reciprocal) infinitely small.  In one of his early
articles,
%
%
Knobloch observes:
\begin{enumerate}\item[]
[Leibniz] distinguished between two infinites, the bounded infinite
straight line, the recta infinita terminata, and the unbounded
infinite straight line, the recta infinita interminata.  He
investigated this distinction in several studies from the year 1676.
Only the first kind of straight lines can be used in mathematics, as
he underlined in his proof of theorem 11 [i.e., Propositio~XI].\, He
assumed a fictive boundary point on a straight halfline which is
infinitely distant from the beginning: a bounded infinite straight
line is a fictitious quantity.  (Knobloch \cite{Kn99}, 1999, p.\;97)
\end{enumerate}
An unbounded infinity, which Leibniz viewed as a contradictory
concept, can be exemplified by the multitude of natural numbers (which
seen as a whole would contradict the part-whole principle; see
Section~\ref{s31b}).  A bounded infinity can be exemplified by an
inassignable (see below) natural number, say~$\mu$, or the multitude
of all natural numbers up to~$\mu$ (and thus \emph{bounded} by~$\mu$).

Leibniz defines an infinitesimal as a ``fraction infiniment petite, ou
dont le denominateur soit un nombre infini'' \cite[p.\;93]{Le02}, such
as~$\frac1\mu$.  Such infinitesimals are routinely used in computing
e.g., the differential ratio~$\frac{dy}{dx}$, as in the passage from
\emph{Cum Prodiisset} (\cite{Le01c}, circ.\;1701), where Leibniz also
mentions the assignable/inassignable dichotomy:
\begin{enumerate}\item[]
[A]lthough we may be content with the assignable quantities
$(d)y$,~$(d)v$,~$(d)z$,
$(d)x$, etc., {\ldots} yet it is plain from what I have said that, at
least in our minds,%
\footnote{For Leibniz on minds, see Section~\ref{s31} and the main
text at note~\ref{f35}.}
the unassignables [\emph{inassigna\-biles} in the original Latin]~$dx$
and~$dy$ may be substituted for them by a method of supposition even
in the case when they are evanescent {\ldots} (Leibniz as translated
in Child~\cite{Ch20}, 1920, p.\,153).
\end{enumerate}
A similar passage appears in \emph{Historia et Origo calculi
differentialis a G. G. Leibnitzio conscripta}.%
\footnote{\;``Although we may be content with the assignable
  quantities~$(d)y$,~$(d)v$,~$(d)z$, and~$(d)x$, since in this way we
  can perceive the whole fruit of our calculus, namely a construction
  using assignable quantities, still it is clear from this that we
  may, at least by feigning, substitute for them the
  unassignables~$dx$,~$dy$ by way of fiction even in the case where
  they vanish, since~$dy:dx$ can always be reduced to~$(d)y:(d)x$, a
  ratio between assignable or undoubtedly real quantities'' (Leibniz
  as translated by RA in \cite[p.\;439]{Ra20}).}
The assignable/inassignable dichotomy was analyzed in the seminal
study by Bos (\cite{Bo74}, 1974).%
\footnote{\label{f30}The dichotomy is also mentioned in Arthur's
  glossary of Leibnizian terms; see main text at note~\ref{f19}.}
See \cite{21a}, Section~4 for a formalisation of the dichotomy in
modern mathematics.

\subsection
{Instances of infinita terminata}
\label{s34b}

In a february 1676 text, Leibniz provided a colorful example of a
bounded infinity as follows:
\begin{enumerate}\item[]
This is a wonder, too: that someone who has lived for infinitely many
years can have begun to live, and that someone who lives for a number
of years that is greater than any finite number can at some time die.
From which it will follow that there is an infinite number.%
\footnote{Later Leibniz will use this example to illustrate the
  impossibility of \emph{infinita terminata} in nature.}
(Leibniz as translated by Arthur in \cite[p.\;51]{Ar01})
\end{enumerate}
In a text \emph{Numeri infiniti} dated 10 april 1676, Leibniz
envisioned the possibility that such numbers - examples of bounded
infinities - may be prime:
\begin{enumerate}\item[]
Two infinite numbers which are not as two finite numbers can be
commensurable, namely, if their greatest common measure [i.e.,
divisor] is a finite number--for instance, if both are prime.
(Leibniz as translated by Arthur in \cite[p.\;87]{Ar01})
\end{enumerate}
The dating of \emph{Numeri infiniti} is significant since Arthur
commits himself to a tight timeframe for an alleged switch in
Leibniz's thinking about infinitesimals.%
\footnote{\label{f32}Thus, Arthur writes: ``For some time, Leibniz
  appears to have hesitated over this interpretation, and as late as
  February 1676 he was still deliberating about whether the success of
  the hypothesis of infinities and the infinitely small in geometry
  spoke to their existence in physical reality too. But by April
  [1676], the syncategorematic interpretation is firmly in place''
  \cite[p.\;559]{Ar13}.  Here Arthur appears to acknowledge that
  Leibniz's hesitation and deliberation concern ``their existence [or
    otherwise] in physical reality''; the ease with which Arthur skips
  from denial of material existence to syncategorematics is
  surprising.  Arthur's dating clashes with Knobloch's scenario
  placing the switch years earlier.
%
}

RA mention the term \emph{terminata} (bounded) three times in their
article \cite{Ra20} and specifically in the context of
\emph{Propositio}\;XI and its \emph{Scholium} in DQA, but do not pay
sufficient attention to the dichotomy of bounded versus unbounded
infinity and fail to appreciate its significance.

Thus, in an analysis of Propositio XI involving the evaluation of a
(finite) area of a region extending to infinity, Leibniz introduces
\begin{enumerate}\item[]
a point~$(\mu)$ at infinitely small distance from the axis.  In this
case, indeed, the straight line~$(\mu)\lambda$ will still be
infinite. (Leibniz as translated by RA in \cite[p.\;421]{Ra20}).  
\end{enumerate}
RA go on to explain this infinity ``in the sense that it \emph{can be
  made} greater than any given quantity (major qualibet assignabili),
yet will be bounded (terminata)'' (ibid.; emphasis added).  However,
there is no source in Leibniz for the clause ``it \emph{can be made}
greater than any given quantity''.  Rather, Leibniz wrote ``it
\emph{is} greater than any given quantity''; RA added the clause ``can
be made greater, etc.''  to help Leibniz conform to the
Ishiguro-syncategorematic interpretation (see Section~\ref{s4}).  RA's
addition distorts Leibniz' intended meaning.

Leibniz concludes the discussion of bounded infinities in the
\emph{Scholium} following \emph{Propositio} XI with the following
comment:
\begin{enumerate}\item[]
D\'eterminer si la nature [\emph{natura rerum}] souffre des
quantit\'es de ce genre est l'affaire du M\'etaphysicien; il suffit au
G\'eom\`etre de d\'emontrer ce qui r\'esulte de leur supposition.%
\footnote{In the original Latin: ``An autem hujusmodi quantitates
ferat natura rerum Metaphysici est disquirere; Geometra sufficit, quid
ex ipsis positis sequatur, demonstrare'' \cite[p.\;98]{Le04b}.
English translation: ``Determining whether nature warrants quantities
of this type is the business of the Metaphysician; for the Geometer it
shall suffice to demonstrate what follows from their supposition.''}
(Leibniz as translated by Parmentier in \cite[p.\,101]{Le04b}).
\end{enumerate}
Here Leibniz observes that, though bounded infinities may not be found
in \emph{rerum natura}, their usefulness to the geometer is
independent of the metaphysician's task of elucidating their relation
to natural phenomena.%
\footnote{\label{f34}Arthur comments as follows on \emph{rerum natura}
  in connection with \emph{Numeri infiniti}: ``In `Infinite Numbers'
  of April 10th any entity such as a line smaller than any assignable,
  or the angle between two such lines, is firmly characterized as
  `fictitious' (A VI, III, 498–99; LLC, 89). There are no such things
  in \emph{rerum natura}, even though they express `real truths'{}''
  (Arthur \cite{Ar09}, 2009, p.\;28).  Arthur goes on to quote Leibniz
  as follows: ``these fictitious entities are excellent abbreviations
  of propositions, and are for this reason extremely useful'' (ibid.).
  Leibniz's comment about infinitesimals being ``excellent
  abbreviations'' is readily appreciated by anyone with experience in
  teaching or research in infinitesimal analysis.  Yet Arthur fails to
  consider the possibility that the absence of infinitesimals in
  \emph{rerum natura} may not imply that they are placemarks for
  quantifier-equipped propositions.  It takes some leap of
  Weierstrassian faith to see Leibniz's comment on abbreviations as
  evidence for an Ishiguroan alternating quantifier reading of
  infinitesimals as logical fictions (see Section~\ref{s4}).
  Leibnizian analogies between infinitesimals and imaginaries
  similarly undercut Ishiguro's reading; see note~\ref{f43}.}
To be useful, they needn't be found in nature:
\begin{enumerate}\item[]
{\ldots} Therefore everything else will also exist in the mind: and in
it everything that I denied to be possible will now be possible.%
\footnote{\label{f35}In the original Latin: ``Erunt ergo in mente et
caetera omnia: et in ea omnia jam fient, fieri quae posse negabam''
\cite[p.\;499]{Le76}.}
(Leibniz, \emph{Numeri infiniti}, as translated by Arthur in
\cite[p.\;91]{Ar01})
\end{enumerate}
The reference is to a discussion on the previous page of a circle as a
fictional polygon with an infinite number of sides:
\begin{enumerate}\item[]
The circle--as a polygon greater than any assignable, as if that were
possible--is a fictive entity, and so are other things of that kind.
(Leibniz as translated by Arthur in \cite[p\;89]{Ar01})
\end{enumerate}

\subsection{Hornangles and inassignables}  
\label{s53}

Leibnizian hornangles exhibit non-Archi\-me\-dean behavior (when
compared to ordinary angles) not easily paraphrasable in
Archi\-me\-dean terms, and shed light on Leibniz's attitude toward
infinitesimals in general.  We provide a geometric illustration of the
phenomenon of magnitudes smaller than all assignables (see
equation~\eqref{e51} in Section~\ref{s42}), in terms of the hornangle
(also known as angle of contact or angle of contingence), much
discussed in 13th--17th century literature.  Thus, Campanus of Novara
(1220--1296) wrote that ``any rectilinear angle is greater than an
infinite number of angles of contingence'' \cite[pp.\;580--581]{Kr}.

The hornangle is the ``angle,'' or crevice, between an arc of a circle
(or a more general curve) at a point~$P$, and its tangent ray~$r$
at~$P$.  According to Leibniz, such a hornangle is smaller than every
ordinary rectilinear angle formed by the ray~$r$ and a secant ray
at~$P$.%
\footnote{\label{f36}In the sense that a sufficiently short subarc
  at~$P$ will lie inside the rectilinear angle.  See Thomason
  (\cite{Th82}, 1982) for a detailed discussion of the non-Archimedean
  nature of hornangles when compared to rectilinear angles.}
Meanwhile, the hornangle is certainly not smaller than all other
hornangles (including itself).

In a 1686 article, Leibniz considered a more general notion of
hornangle, or angle of contact, between a pair of curves with a common
tangent at~$P$.%
\footnote{From the viewpoint of modern geometry the angle between a
  pair of tangent curves is zero by definition, but Leibniz envisioned
  a more general notion of angle.}
Leibniz defined the osculating circle of a curve at~$P$ as the circle
forming the least angle of contact with the curve:
\begin{enumerate}\item[]
Circulus autem ille lineam propositam ejusdem plani in puncto
proposito osculari a me dicitur, qui minimum cum ea facit angulum
contactus.   (Leibniz \cite{Le86}, 1686)
\end{enumerate}
Choosing the least one among the angles of contact, as Leibniz does
here, presupposes that such ``angles'' are nonzero.  A few lines
later, Leibniz renames such entities ``angle of osculation'' to
distinguish them from angles of contact in the classical sense
(crevice between curve and tangent ray).

Leibniz explicitly compares angles of contact to inassignables in the
1671 text \emph{Theoria motus abstracti} (TMA).  In this early text,
Leibniz still composes the continuum of indivisible infinitesimals
which he refers to as ``points.''  Speaking of the space filled by an
infinitesimal motion of a body, Leibniz writes:
\begin{enumerate}\item[]
{\ldots} this space is still inassignable {\ldots}\;although the ratio
of a point of a body {\ldots}\;to the [part]%
\footnote{The term occurring in the original, \emph{point}, seems to
  be a misprint.}
of space it fills when moving, is as an angle of contact to a
rectilinear angle {\ldots} (Leibniz as translated by Loemker
\cite[p.\,140]{Le89}, TMA, Predemonstrable Foundation 13, A VI.ii N41)
\end{enumerate}
In his mature period, Leibniz's view of the continuum changes (see
Section~\ref{anto2}), but the non-Archimedean behavior of hornangles
in comparison with ordinary rectilinear/assignable angles remains.
There are at least three indications that Leibniz viewed hornangles as
being nonzero.  

\medskip
\textbf{1.}  Choosing the smallest one of the angles of contact as
Leibniz does in 1686 presupposes that they are nonzero.  

\medskip
\textbf{2.}  In a 1696 letter, Leibniz makes it clear that angles of
contact are nonzero:
\begin{enumerate}\item[]
Our infinitesimal calculus allows us to see that one can only ignore
\emph{differentias incomparabiliter minores rebus differentiatis} (the
differences which are incomparably smaller than the differentiated
things).  So, it does not follow that there is no considerable
difference between the degrees of force in the object from each blow
by gravitating matter.  Otherwise, this would be as if one wanted to
prove that angles of contact do not differ amongst themselves because
they do not compare with rectilinear angles.  (Leibniz to Papin, March
1696, A III, 6, 698).
\end{enumerate}
Here Leibniz uses the term \emph{incomparable} in the technical sense
(akin to \emph{inassignable}) used in his formulation of the violation
of the Archimedean property in a 1695 letter to L'H\^opital.%
\footnote{See the quotation in note~\ref{f71}.}
Of course, they can still be \emph{compared} to rectilinear angles.%
\footnote{As explained in note~\ref{f36}.}

\medskip
\textbf{3.}  A careful analysis of the two passages in (\cite{Le86},
1686) where Leibniz asserts that they are ``nothing'' reveal that his
intention is to say that they are negligible vis-\`a-vis incomparably
larger quantities: first, angles of contact are negligible vis-\`a-vis
rectilinear angles, and second, angles of osculation are negligible
vis-\`a-vis angles of contact.

\medskip
We note the unavailability of the option of representing a hornangle
by rectilinear angles -- either arranged in a sequence or assorted
with logical quantifiers.  The non-Archimedean behavior of Leibnizian
hornangles is not easily paraphrasable in Archimedean terms,
suggesting similar behavior of Leibnizian inassignable infinitesimals.
A similar situation exists with respect to comparison of
infinitesimals and imaginary roots.%
\footnote{See notes \ref{f34}, \ref{f43}, and \ref{f50}.}

Infinitesimals, hornangles, and imaginaries are well-founded fictions
that facilitate the art of discovery.  Leibnizian well-founded
fictions are at most accidentally impossible (they are not
contradictory; see Section~\ref{s47}).  Meanwhile, there are entities
that Leibniz sometimes refers to as fictions, such as infinite wholes,
which are absolutely impossible (contradictory).%
\footnote{\label{f42}Thus, analyzing the area under the hyperbola,
  Leibniz concludes: ``By this argument it is concluded that the
  infinite is not a whole, but only a fiction, since otherwise the
  part would be equal to the whole'' (A VII 3, 468; october 1674).
  Arthur \cite[p.\;557]{Ar13} and Rabouin--Arthur \cite[p.\;405]{Ra20}
  quote this passage but fail to account for the fact that Leibniz
  never refers to such entities as \emph{well-founded} fictions.  For
  more details see \cite{22b}.}
This crucial distinction is not sufficiently taken into account by
advocates of the Ishiguro-syncategorematic interpretation.
%
%

\section{Well-founded fictions: the philosophy}
\label{s2c}

\subsection{Instantiation in \emph{rerum natura}?}
\label{s31}

On a number of occasions, Leibniz spoke of infinitesimals as not
necessarily \emph{found in nature}.  Thus, in the \emph{Specimen
  Dynamicum} of 1695, he concluded his analysis of infinite degrees of
impetus as follows:
\begin{enumerate}\item[]
Hence the nisus is obviously twofold, an elementary or infinitely
small one which I also call a \emph{solicitation} and one formed by
the continuation or repetition of these elementary impulsions, that
is, the impetus itself.  But I do not mean that these mathematical
entities \emph{are really found in nature} as such but merely that
they are means of making accurate calculations of an \emph{abstract
  mental kind}.  (Leibniz \cite{Le95d} as translated by Loemker in
\cite[p.\;438]{Le89}; emphasis on \emph{solicitation} in the original;
emphasis on \emph{really found in nature} and \emph{abstract mental
  kind} added)
\end{enumerate}
Over a decade later, des Bosses questioned Leibniz concerning the
above passage from \emph{Specimen Dynamicum}.  Leibniz responded as
follows concerning the fictionality of infinitesimals:
\begin{enumerate}\item[]
For I consider both [infinitely small and infinitely large magnitudes]
as \emph{fictions of the mind}, due to abbreviated ways of speaking,
which are suitable for calculation, in the way that imaginary roots in
algebra are.%
\footnote{In the original Latin: ``Utrasque enim per modum loquendi
compendiosum pro mentis fictionibus habeo, ad calculum aptis, quales
etiam sunt radices imaginariae in Algebra'' \cite[p.\;32]{Le06}.  
%
%

Schubring quotes this sentence and claims that ``Leibniz stressed that
for himself the infinitely small quantities were not really existing
mathematical quantities, but only \emph{fictions} that had their uses
in the course of calculus'' (\cite[p.\,171]{Sc05}; emphasis in the
original).  It is correct that Leibniz referred to infinitesimals as
\emph{fictions}, but Schubring's claim that to Leibniz they ``were not
really existing mathematical quantities'' is unsupported and at best
ambiguous.  For a detailed critique of Schubring's book see
\cite{17e}.}
(Leibniz as translated by Look and Rutherford in \cite[p.\;33]{Le06};
emphasis added)
\end{enumerate}
Arthur quotes this passage in (\cite{Ar18}, 2018, p.\,176) and
concludes: ``This is the syncategorematic infinitesimal described
above'' (ibid.).%
\footnote{Similar claims concerning this Leibnizian passage appeared
five years earlier in Arthur (\cite{Ar13}, 2013, p.\;555).}
Is it indeed?  First, in this passage Leibniz refers to \emph{both}
infinitesimals and imaginary roots as \emph{compendia}
(abbreviations), undermining an Ishiguro-syncatego\-re\-matic reading
(see Section~\ref{s4}) since it is unavailable for imaginary roots.%
\footnote{\label{f43}Since Leibniz describes both as abbreviations,
  the absence of any plausible alternating quantifier account of
  imaginary roots in terms of more ordinary quantities suggests that
  it was not Leibniz's intention in the case of infinitesimals,
  either.  A similar comparison occurs in a 2\;february 1702 letter to
  Varignon; see main text at note~\ref{f50}.}
Furthermore, des Bosses, in formulating his question on 2 march 1706,
specifically raised the possibility of the syncategorematic infinite:
\begin{enumerate}\item[]
I would have conjectured that the infinite that you add can be
confined to the syncategorematic; {\ldots} (des Bosses as translated
by Look and Rutherford in \cite[p.\;27]{Le07})
\end{enumerate}
In his answer nine days later, on 11 march 1706, Leibniz says not a
word about the syncategorematic infinite, and rather speaks of
fictions of the mind as quoted above.

Levey reproduces a longer passage from the 11\;march 1706 letter
containing the one we quoted, and claims that it supports the
syncategorematic reading (\cite{Le21}, 2021, p.\,146).  However, such
a claim overlooks the fact that Leibniz specifically ignored des
Bosses' question about the syncategorematic infinite, as noted above.
Levey concludes his Section\;2.1 on infinitesimals as follows:
\begin{enumerate}\item[]
They can be replaced by proofs given in Leibniz's updated style of
Archimedes if full rigor is wanted, and the mathematics in which they
figure is not committed to the existence of `actual' infinitesimals
\emph{in\;nature} for its justification.  (op.\;cit., p.\,148;
emphasis added)
\end{enumerate}
As Levey appears to acknowledge in his concluding sentence, it is
\emph{in\;nature} that there may be no infinitesimals; Leibniz viewed
them as \emph{mentis fictiones}.  Levey's stated conclusion is mainly
in accord with our reading of Leibnizian infinitesimals.  Meanwhile,
the Ishiguro-syncategorematic hypothesis pursued elsewhere in Levey's
Section\;2.1 remains unsupported.

\subsection
{Well-founded fictions in relation to \emph{rerum natura}}
\label{s33}

The fictional nature of infinitesimals is a constant theme in
Leibnizian thought.  Thus, in a 14\;april\;1702 letter to Varignon,
Leibniz wrote:
\begin{enumerate}\item[]
[L]es infinis et infiniment petits pourroient estre pris pour des
fictions, semblables aux racines imaginaires, sans que cela d\^ut
faire tort \`a nostre calcul, ces fictions estant utiles et fond\'ees
en realit\'e.%
\footnote{Translation: ``Infinities and the infinitely small could be
taken for fictions, similar to imaginary roots, without it causing
harm to our calculus, these fictions being useful and founded in
reality.''}
\cite[p.\;98]{Le02b}
\end{enumerate}
In this passage, Leibniz asserts that interpreting the infinite and
infinitely small as useful fictions would cause no harm to `our
calculus'.%
\footnote{Inexplicably, Schubring's rendition attaches the opposite
  meaning to the passage: ``Infinities and infinitely small quantities
  could be taken as fictions, similar to imaginary roots, except that
  \emph{it would make our calculations wrong}, these fictions being
  useful and based in reality'' (Leibniz as rendered by Schubring in
  \cite[p.\,171]{Sc05}; emphasis added).  Only Schubring's rendition
  makes the calculations wrong.}
In a 20\;june 1702 letter to Varignon, Leibniz wrote:
\begin{enumerate}\item[]
Entre nous je crois que Mons.\;de Fontenelle
%
%
{\dots}\ en a voulu railler, lorsqu'il a dit qu'il vouloit faire des
elemens metaphy\-siques de nostre calcul.  Pour dire le vray, je ne
suis pas trop persuad\'e moy m\^eme, qu'il faut considerer nos infinis
et infiniment petits autrement que comme \emph{des choses ideales ou
  comme des fictions bien fond\'ees}.%
\footnote{Translation: ``Between us, I believe that Mr.\;Fontenelle
{\ldots} was joking when he said that he wished to develop
metaphysical elements of our calculus.  To tell the truth, I am not
myself persuaded that it is necessary to consider our infinities and
the infinitely small as something other than ideal things or
well-founded fictions.''}
(Leibniz \cite{Le02c}, 1702, p.\,110; emphasis added)
\end{enumerate}
In this case Leibniz adds a qualification: infinitesimals are
\emph{well-founded} fictions (`fictions bien fond\'ees').  It is
difficult to see how fictions that are, according to Leibniz,
\emph{well-founded} could also be, as per RA, contradictory.%
\footnote{\label{f48}RA attempt to sidestep the difficulty by claiming
  that contradictory infinitesimals can be used ``under certain
  conditions'' but don't specify the latter; see main text at
  note~\ref{f24}.  In his letter to des Bosses dated 11\;march\;1706,
  Leibniz describes the rainbow as a well-founded phenomenon
  \cite[p.\;35]{Le06}.  Since Leibniz describes infinitesimals using
  the same expression ``well-founded,'' the question arises whether RA
  would be prepared to claim that to Leibniz, rainbows similarly were
  contradictory notions that can be used under certain unspecified
  conditions.  Occasionally, Leibniz refers to infinite wholes as
  ``fictions,'' as in a 1674 text, but he never refers to them as
  \emph{well-founded} fictions; see note~\ref{f42}.}
On the contrary, for Leibniz consistency is a requirement for a
well-founded fiction and thus for mathematical existence.  In this
respect Leibniz is closer than many mathematicians of subsequent
generations to Hilbert’s formalism (where existence depends on
consistency alone).  By ascribing to Leibniz the use of contradictory
concepts, RA rule out an interpretation whereby Leibniz, like Hilbert,
views conception of mathematical existence as consistency; see
\cite{13f} for details.

Leibniz expressed similar sentiments concerning ideal notions in a
2\;february\;1702 letter to Varignon:%
\footnote{\label{f49}RA's interpretation of the letter is analyzed in
Section~\ref{s54}.}
\begin{enumerate}\item[]
D'o\`u il s'ensuit, que si quelcun n'admet point des lignes infinies
et infiniment petites \`a la rigueur metaphysique et comme des choses
reelles, il peut s'en servir seurement comme des \emph{notions
ideales} qui abregent le raisonnement,%
\footnote{\label{f50} ``Notions ideales qui abregent le raisonnement''
  is Leibniz's French equivalent of the Latin \emph{compendia}; see
  note~\ref{f43}.  In this passage also, infinitesimals and
  imaginaries are both described as abbreviations that facilitate
  reasoning.}
semblables \`a ce qu'on appelle racines imaginaires dans l'analyse
commune (comme par exemple~$\sqrt{-2}$), {\ldots}%
\footnote{Loemker's translation: ``It follows from this that even if
someone refuses to admit infinite and infinitesimal lines in a
rigorous metaphysical sense and as real things, he can still use them
with confidence as ideal concepts which shorten his reasoning, similar
to what we call imaginary roots in the ordinary algebra, for example,
$\sqrt{-2}$\,'' \cite[p.\;543]{Le89}.}
(\cite[p.\;92]{Le02}; emphasis added)
\end{enumerate}
The idea of imagining infinitesimals also appears in a 30\;march\;1699
letter to Wallis, where Leibniz rejects Wallis' position that
infinitesimals are \emph{nothings}:
\begin{enumerate}\item[]
[F]or the calculus it is useful \emph{to imagine}%
\footnote{Possibly: \emph{to feign}; in the original Latin:
\emph{fingere}.  See also main text at note~\ref{f35}.}
infinitely small quantities, or, as Nicolaus Mercator called them,
infinitesimals, such that when at least the assignable ratio between
them is sought, they precisely may not be taken to be nothings.
(Leibniz as translated by Beeley in \cite[note\;38]{Be08}; emphasis
added)
\end{enumerate}
In the same letter, Leibniz makes a revealing comment concerning the
status of inassignables:
\begin{enumerate}\item[]
Whether inassignable quantities are real or fictions, I will not argue
for now; it is enough that they serve as a help for thinking, and that
they always carry a proof with them, with only the style changed; and
so I have noted, that if anyone should substitute incomparably or
sufficiently small (quantities) for infinitely small ones, I do not
object.  (Leibniz to Wallis, 30 march 1699, GM \cite{Ge50} IV, 63;
translation ours)
\end{enumerate}
Here Leibniz refuses to commit himself to either a realist or a
fictionalist position.  Beeley offers the following intriguing
speculation:
\begin{enumerate}\item[]
This of course opens up the whole question of whether Leibniz really
held that infinitesimals could exist in nature.  On some occasions he
does indeed seem to be denying their existence.  But I think that we
need to be careful here, because denial of the existence of
infinitesimals is generally coupled with the argument that the success
of the calculus does not depend on metaphysical discussions concerning
reality.  When he makes such claims, this seems to be no more than a
get-out clause vis-\`a-vis opponents who seek to provide metaphysical
arguments against his calculus.  Seen within the context of Leibniz'
dynamics, particularly in respect of dead force (vis mortua) it is
evident that he must be committed in some way to the existence of
infinitesimals.%
\footnote{For related comments on dead force by Garber see
  Section~\ref{anto}.}
(Beeley \cite{Be15}, 2015, p.\;42)
\end{enumerate}
Our main arguments in the present text are independent of resolving
Beeley's \mbox{`get-out clause'} hypothesis.  On many occasions,
Leibniz did reject infinitesimal creatures.  Thus, in a 20 june 1702
letter to Varignon, Leibniz wrote:
\begin{enumerate}\item[]
Je croy qu'il n'y a point de creature au dessous de la quelle il n'y
ait une infinit\'e de creatures, cependant je ne crois point qu'il y
en ait, ny m\^eme qu'il y en puisse avoir d'infiniment petites et
c'est ce que je crois pouvoir demonstrer.
%
%
Il est que les substances simples (c'est \`a dire qui ne sont pas des
estres par aggregation) sont veritablement indivisibles, mais elles
sont immaterielles, et ne sont que principes d'action.  (Leibniz
\cite{Le02c}, 1702, p.\,110)
\end{enumerate}
In the same vein, in an 11 march 1706 letter to des Bosses, Leibniz
wrote:
\begin{enumerate}\item[]
Yet you see that it should not be concluded from this that an
infinitely small portion of matter (such as does not exist) must be
assigned to any entelechy, even if we usualy rush to such conclusions
by a leap.  (Leibniz as translated by Look and Rutherford in
\cite[p.\;35]{Le06})
\end{enumerate}
%

\subsection
{Theory of knowledge: two types of impossibility}
\label{s47}

The terms ``contradictory'' and ``impossible'' have different meanings
for Leibniz.  In Leibnizian theory of knowledge, the fact that
something is (1)\,\emph{not possible} does not mean that it is
(2)\,\emph{absolutely impossible or contradictory}.

Leibniz introduces a related distinction on several occasions.  Thus,
in the \emph{Confessio Philosophi} (Leibniz \cite{Le72}, 1672/3,
p.\,128), he refers to (1) as ``impossible by accident'' and contrasts
it with (2)\;``abso\-lute impossibility'' i.e., contradiction.  He
gives the examples of a species with an odd number of feet, and an
immortal mindless creature, which are, according to him,
\emph{harmoniae rerum adversa} i.e., ``contrary to the harmony of
things'' (trans.\;Sleigh in \cite[p.\;57]{Le04c}), but \emph{not}
contradictory.

In his 1683 text \emph{Elementa nova matheseos universalis}, Leibniz
explains that some mathematical operations cannot be performed in
actuality, but nonetheless one can exhibit ``a construction in our
characters'' (\emph{in nostris characteribus} \cite[p.\;520]{Le83}) --
meaning that one can carry out a formal calculation, such as those
with imaginary roots, regardless of whether the mathematical notions
involved idealize anything in nature.  Leibniz goes on to discuss in
detail the cases of imaginary roots and infinitesimals.  For
convenience of reference, we labeled four passages [A], [B], [C], [D].
Leibniz mentions infinitesimals in passage~[C]:
\begin{enumerate}\item[]
[A] And some extractions of roots are such that roots are surd and
they do not exist in \emph{natura rerum}, and we call them imaginary,
and the problem is impossible, as when our analysis shows that the
requested point must be exhibited by the intersection of a specific
circle and a specific straight line, in which case it may happen that
this circle by no means reaches this line, and then the intersection
is imaginary {\ldots}

\medskip\noindent [B] There is a big difference between imaginary
quantities, or those impossible by accident, and absolutely impossible
ones, which involve a contradiction: e.g., when it is found that
solving a problem requires that 3 be equal to 4, which is absurd.

\medskip\noindent
[C] But imaginary quantities, or [quantities] impossible by accident,
namely quantities that cannot be exhibited for lack of a sufficient
condition, which is required for having an intersection, can be
compared with infinite and infinitely small quantities, which are
generated in the same way.  {\ldots}

\medskip\noindent [D] And it is true that calculus necessary leads to
them, and people who are not sufficiently expert in such matters get
entangled and think they have reached an absurdity [\emph{absurdum}].
Experts know instead that this apparent impossibility
[\emph{apparentem illam impossibilitatem}] only means that a parallel
line is traced instead of a straight line making the required angle,
and this parallelism is the required angle, or quasi-angle.%
\footnote{\;``Et quaedam extractiones tales sunt, ut radices illae
  surdae nec in natura rerum extent, tunc dicuntur imaginariae, et
  problema est impossibile, ut cum analysis ostendit punctum quaesitum
  debere exhiberi per intersectionem certi circuli et certae rectae,
  ubi fieri potest ut ille circulus ad illam rectam nullo modo
  perveniat, et tunc intersectio erit imaginaria {\ldots} Multum autem
  interest inter quantitates imaginarias, seu impossibiles per
  accidens, et impossibiles absolute quae involvunt contradictionem,
  ut cum invenitur ad problema solvendum opus esse, ut fit 3 aequ.\;4
  quod est absurdum. Imaginariae vero seu per accidens impossibiles,
  quae scilicet non possunt exhiberi ob defectum sufficientis
  constitutionis ad intersectionem necessariae, possunt comparari cum
  Quantitatibus infinitis et infinite parvis, quae eodem modo
  oriuntur.  {\ldots} quod adeo verum est ut saepe calculus ad eas
  necessario ducat, ubi harum rerum nondum satis periti mire
  torquentur et in absurdum se incidisse putant.  Intelligentes vero
  sciunt apparentem illam impossibilitatem tantum significare, ut loco
  rectae angulum quaesitum facientis ducatur parallela; hunc
  parallelismum esse angulum illum seu quasi angulum quaesitum.''}
(Leibniz \cite{Le83}, 1683, pp.\;520--521)%
\footnote{This passage was translated into French by Rabouin in
  \cite[pp.\,107--109]{Le18}.}
\end{enumerate}
Passages [B] and [C] indicate that (natural instantiations of)
infinitesimals are only \emph{impossible by accident}, or \emph{only
  apparently impossible}.  Therefore there is no grounds for
attributing absolute impossibility or contradiction to them, or for
lumping them with infinite wholes, as per RA (see Section~\ref{s43}).
Infinite wholes are absolutely impossible because they are contrary to
the part-whole principle (see Section~\ref{s31b}), which is a
necessary truth.
%
%
Thus, in a 1678 letter to Elisabeth \cite{Le78}, Leibniz described the
concept of ``the number of all possible units''%
\footnote{Here ``number'' refers to cardinality as per Ishiguro; see
  note~\ref{f21}.}
as impossible \cite[p.\;238]{Le89a}.  In his \emph{Historia et Origo}
\cite{Ge46}, he presents a derivation of the part-whole principle from
the principle of identities and the definitions of whole and part:
\begin{enumerate}\item[]
[T]hat mighty axiom, `The whole is greater than its part', could be
proved by a syllogism of which the major term was a definition and the
minor term an identity.  (Leibniz as translated by Child in
\cite[pp.\;29--30]{Ch20})
\end{enumerate}
RA seek to undercut infinitesimals on the grounds that their
definition, in terms of violation of Euclid's Definition V.4,%
\footnote{See note~\ref{f71}.}
is only nominal, and claim the following:
\begin{enumerate}\item[]
Although this concept contains a contradiction, other subsidiary
concepts contained in it may permit the deri\-vation of true
entailments.  \cite[p.\;406]{Ra20}
\end{enumerate}
Does the concept of infinitesimal contain a contradiction as claimed
by RA?  In his 1686 comments that shed light on the
Leibnizian theory of knowledge, Leibniz wrote merely that a nominal
definition \emph{could} harbor contradictions, not that it must do so:
\begin{enumerate}\item[]
Et tandis qu'on n'a qu'une definition nominale, on ne se s\c{c}auroit
asseurer des consequences qu'on en tire; car si elle cachoit quelque
contradiction ou impossibilit\'e, on en pourroit tirer des conclusions
oppos\'ees. (Leibniz \cite{Le86b}, 1686, pp.\,1568--1569)
\end{enumerate}
%
%
With regard to both imaginaries and infinitesimals, Leibniz makes it
clear in the 1683 passage cited above (paragraph [C]) that their
natural instantiations are (at most) only conditionally or
accidentally impossible, rather than contradictory.  Esquisabel and
Raffo\;Quintana examine the issue and reach the following conclusion:
\begin{quote}
[U]nlike the infinite number or the number of all numbers, for Leibniz
infinitary concepts do not imply any contradiction, although they may
imply paradoxical consequences.  \cite[p.\;641]{Es21}%
\footnote{\label{f60b}Esquisabel and Raffo Quintana clarify: ``[W]e
  disagree with the reasons [Rabouin and Arthur] gave for the
  Leibnizian rejection of the existence of infinitesimals, and in our
  opinion the texts they refer to in order to support their
  interpretation are not convincing.  Since we argue that Leibniz did
  not consider the concept of infinitesimal as
  \emph{self-contradictory}, we try to provide an alternative
  conception of impossibility''\;\cite[p.\;620]{Es21}.}
\end{quote}
If ``the principle of non-contradiction [is] the principle of the
minimal condition of intelligibility'' (Grosholz and
Yakira\;\cite{Gr98}, 1998, p.\;44) then one can easily perceive why
infinitesimals, unlike infinite wholes, are ubiquitous in Leibniz's
mathematical oeuvre.

\subsection
{Infinitesimals, infinite wholes, and nominal definitions} 

To summarize, in a 1672/3 text \emph{Confessio Philosophi}, Leibniz
speaks of the distinction between accidental impossibility and
absolute impossibility (equivalent to contradiction), and mentions two
examples of accidental impossibility, specifying that they are
contrary to the principle of the \emph{harmoniae rerum}.

Furthermore, in a 1683 text \emph{Elementa nova matheseos universalis}
Leibniz mentions that imaginary roots do not exist in \emph{rerum
  natura}, and contrasts imaginary quantities impossible by accident,
on the one hand, and absolutely impossible entities involving a
contradiction such as~$3=4$, on the other.  Leibniz goes on to compare
imaginary quantities (impossible by accident) to infinitely small
quantities, and points out that people who are not sufficiently expert
tend to confuse apparent impossibility with absurdity.

Finally, in a 1686 text \emph{Discours de Metaphysique} Leibniz speaks
of nominal definitions and warns that they \emph{might} harbor
contradictions (but not that they must necessarily do so).  Even if
the definition of an infinitesimal as violating Euclid's
Definition\;V.4 were nominal (as RA claim), it would follow at most
that infinitesimals, possibly contrary to the principle of the
\emph{harmoniae rerum}, may not idealize anything in nature, a state of
affairs that can be described as accidental impossibility.  Meanwhile,
infinite wholes are contrary to the part-whole principle (which
Leibniz consistently takes to be a necessary truth), and therefore do
involve an absurdity.

\section{Conclusion}
\label{s8}

Bassler (\cite{Ba98}, 1998), Arthur (\cite{Ar19b}, 2019) and
Rabouin--Arthur (\cite{Ra20}, 2020) attempt to enlist Leibniz's 1716
letter to Masson in support of an Ishiguro-syncategorematic reading of
the Leibnizian calculus, claiming that its procedures are in keeping
with the Archimedean axiom.%
\footnote{Meanwhile, Rabouin solo \cite{Ra15b} expressed a more
cautious position with regard to the 1716 letter; see
Section~\ref{s54}, item 3.}
Rabouin and Arthur assimilate Leibnizian infinitesimals to infinite
wholes and surmise that Leibniz viewed infinitesimals as
contradictory.  Such a reading rules out an interpretation whereby
Leibniz viewed mathematical existence as consistency, as did Hilbert.

Contextualizing the Leibnizian comments on the calculus in the 1716
letter suggests a different reading.  The letter is consistent with
Leibniz's position in 1695 \cite{Le95a}, \cite{Le95b} that
infinitesimals violate Euclid's Definition\;V.4 when compared to
assignable quantities (hardly in keeping with the Archimedean axiom).%
\footnote{\label{f71}Thus, Breger writes: ``In a letter to L'H\^opital
  of 1695, Leibniz gives an explicit definition of incomparable
  magnitudes: two magnitudes are called incomparable if the one cannot
  exceed the other by means of multiplication with an arbitrary
  (finite) number, and he expressly points to Definition\;5 of the
  fifth book of Euclid'' \cite[pp.\;73--74]{Br17}.}
Like imaginary roots, infinitesimals are useful fictions that are at
most accidental impossibilities that violate the principle of the
harmony of things if they do not idealize anything in nature.
Infinitesimals are therefore to be contrasted with infinite wholes,
which are absolute impossibilities since they contradict the
part-whole principle, a necessary truth in Leibnizian thought.  Our
conclusions are compatible with those of Esquisabel and Raffo
Quintana, who similarly reject the contention by Rabouin and Arthur
that Leibniz viewed infinitesimals as contradictory.

Leibniz wrote to Huygens as follows:
\begin{enumerate}\item[]
[Nieuwentijt] me fait une objection sur un point qui m'est commun avec
Messieurs Fermat, Barrow, Newton et tous les autres, qui ont
raisonn\'e sur les grandeurs infiniment petites.  (Leibniz
\cite{Le95c}, 1695)
\end{enumerate}
Here Leibniz asserts the existence of a point in common between his
reasoning with infinitesimals and that of his illustrious
predecessors.  To maintain the syncategorematic hypothesis, its
proponents Arthur, Bassler, Levey, and Rabouin would face an
unenviable alternative: either (1) claim that Leibniz was untruthful
in his letter to Huygens, or (2)\;argue that, in Leibniz's view, the
infinitesimals of Fermat, Barrow, and Newton were similarly
syncategorematic.

The 1716 letter to Masson provides no basis for doubting that Leibniz
based his calculus on a special kind of small magnitudes that he
viewed as \emph{mentis fictiones}.


\section*{Acknowledgments}   

We are grateful to Roger Ariew and Philip Beeley for helpful
suggestions, to David Schaps for expert advice on Greek spelling and
Latin usage, and to Stephan Meier-Oeser for providing access to the
autograph manuscripts of the letter from Leibniz to Masson relied upon
in note~\ref{f6}.


\begin{thebibliography}{AII}


\bibitem{An15} Antognazza, M.\, The hypercategorematic infinite.
  \emph{The Leibniz Review} \textbf{25} (2015), 5--30.

\bibitem{Ar01} Arthur, R. (Tr.) \emph{The Labyrinth of the Continuum.
  Writings on the Continuum Problem, 1672--1686.  G. W. Leibniz.}  New
  Haven: Yale University Press, 2001.


\bibitem{Ar09} Arthur, R.\, Actual infinitesimals in Leibniz's early
  thought, 11--28 in \emph{The Philosophy of the Young Leibniz},
  Studia Leibnitiana Sonderhefte 35, ed. Mark Kulstad, Mogens Laerke
  and David Snyder, Stuttgart: Franz Steiner, 2009.



\bibitem{Ar13} Arthur, R.\, Leibniz's syncategorematic infinitesimals.
\emph{Archive for History of Exact Sciences} \textbf{67} (2013),
no.\;5, 553--593.


\bibitem{Ar18} Arthur, R.\, Leibniz's syncategorematic actual
  infinite.  In \emph{Infinity in Early Modern Philosophy}, Ohad
  Nachtomy and Reed Winegar (Eds.), Cham: Springer, 2018,
  pp.\,155--179.

\bibitem{Ar19b} Arthur, R.\, Leibniz in Cantor's Paradise.\, In
  \emph{Leibniz and the Structure of Sciences.  Modern Perspectives on
    the History of Logic, Mathematics, Epistemology}, 71--109.
  Ed.\;V. De Risi, Boston Studies in Philosophy and History of
  Science, Cham: Springer, 2019.

\bibitem{Ar19c} Arthur, R.\, Review of \emph{Varieties of
  continua--from regions to points and back} by Hellman and Shapiro
  \cite{He18}.  \emph{Philos. Math. (3)} \textbf{27} (2019), no.\,1,
  148--152.



\bibitem{17b} Bair, J.; B{\l}aszczyk, P.; Ely, R.; Henry, V.; Kanovei,
  V.; Katz, K.; Katz, M.; Kutateladze, S.; McGaffey, T.; Reeder, P.;
  Schaps, D.; Sherry, D.; Shnider,\;S.\, Interpreting the
  infinitesimal mathematics of Leibniz and Euler.  \emph{Journal for
    General Philosophy of Science} \textbf{48} (2017), no.\;2,
  195--238.  See \url{http://dx.doi.org/10.1007/s10838-016-9334-z} and
  \url{https://arxiv.org/abs/1605.00455}


\bibitem{21a} Bair, J.; B{\l}aszczyk, P.; Ely, R.; Katz, M.;
  Kuhlemann, K.\, Procedures of Leibnizian infinitesimal calculus: An
  account in three modern frameworks.  \emph{British Journal for the
    History of Mathematics} \textbf{36} (2021).
  \url{https://doi.org/10.1080/26375451.2020.1851120},
  \url{https://arxiv.org/abs/2011.12628}



\bibitem{16a} Bascelli, T.; B{\l}aszczyk, P.; Kanovei, V.; Katz, K.;
  Katz, M.; Schaps, D.; Sherry, D.\, Leibniz versus Ishiguro: Closing
  a Quarter Century of Syncategoremania.  \emph{HOPOS: The Journal of
    the International Society for the History of Philosophy of
    Science} \textbf{6} (2016), no.\,1, 117--147.
  \url{http://doi.org/10.1086/685645},
  \url{https://arxiv.org/abs/1603.07209}


\bibitem{Ba98} Bassler, O.\, Leibniz on the indefinite as infinite.
\emph{The Review of Metaphysics} \textbf{51} (1998), no.\;4, 849--874.



\bibitem{Be08} Beeley, P.\, Infinity, infinitesimals, and the reform
of Cavalieri: John Wallis and his critics.  In Goldenbaum--Jesseph
\cite{Go08}, pp.\;31--52, 2008.


\bibitem{Be15} Beeley, P.\, Leibniz, philosopher mathematician and
mathematical philosopher.  In \emph{G. W. Leibniz, Interrelations
between mathematics and philosophy}, pp.\;23--48, Archimedes, 41,
Dordrecht: Springer, 2015.


\bibitem{17e} B{\l}aszczyk, P.; Kanovei, V.; Katz, M.; Sherry, D.\,
  Controversies in the foundations of analysis: Comments on
  Schubring's Conflicts.  \emph{Foundations of Science} \textbf{22}
  (2017), no.\,1, 125--140.
  \url{http://doi.org/10.1007/s10699-015-9473-4},
  \url{https://arxiv.org/abs/1601.00059}


\bibitem{Bo74} Bos, H.\, Differentials, higher-order differentials and
the derivative in the Leibnizian calculus.  \emph{Archive for History
of Exact Sciences} \textbf{14} (1974), 1--90.



\bibitem{Bo91} Bosinelli, F.\, \"Uber Leibniz' Unendlichkeitstheorie.
\emph{Studia Leibnitiana} \textbf{23} H.\;2 (1991), 151--169.




\bibitem{Br16} Breger, H.\, Das Kontinuum bei Leibniz.  In W. Li
  (Ed.), \emph{Kontinuum, analysis, informales---Beitr\"age zur
  Mathemathik und Philosophie von Leibniz}, pp.\,115--126.
  Berlin--Heidelberg: Springer, 2016.



\bibitem{Br17} Breger, H.\, On the grain of sand and heaven's
  infinity.  In \emph{`F\"ur unser Gl\"uck oder das Gl\"uck anderer'
    Vortr\"age des X. Intemationalen Leibniz-Kongresses Hannover,
    18.-23. Juli 2016}, Wenchao Li (ed.), in collaboration with Ute
  Beckmann, Sven Erdner, Esther-Maria Errulat, J\"urgen Herbst, Helena
  Iwasinski und Simona Noreik, Band VI, Hildesheim--Zurich--New York:
  Georg Olms Verlag, 2017, pp.\;64--79.



\bibitem{Br99} Brown, S.\,
%
%
Two papers by John Toland. \emph{I Castelli di Yale} \textbf{IV (4)}
(1999), 55--79.





\bibitem{Ch20} Child, J. (Ed.)  \emph{The early mathematical
  manuscripts of Leibniz}.  Translated from the Latin texts published
  by Carl Immanuel Gerhardt with critical and historical notes by
  J. M. Child.  Chicago--London: The Open Court Publishing, 1920.
  Reprinted by Dover in 2005.




\bibitem{Eh06} Ehrlich, P.\, The rise of non-Archimedean mathematics
  and the roots of a mis\-conception. I. The emergence of
  non-Archimedean systems of magnitudes.  \emph{Archive for History of
    Exact Sciences} \textbf{60} (2006), no.\,1, 1--121.



\bibitem{Es21} Esquisabel, O.; Raffo Quintana, F.\, Fiction,
  possibility and impossibility: Three kinds of mathematical fictions
  in Leibniz’s work.  \emph{Archive for History of Exact Sciences}
  \textbf{75} (2021), no.\;6, 613--647.
  \url{https://doi.org/10.1007/s00407-021-00277-0}


\bibitem{Ga08} Garber, D.\, Dead force, infinitesimals, and the
  mathematicization of nature.  In Goldenbaum--Jesseph \cite{Go08},
  2008, pp.\;281--306.


\bibitem{Ge46} Gerhardt, C. (ed.)  \emph{Historia et Origo calculi
  differentialis a G. G. Leibnitio conscripta}.  Hannover: Hahn, 1846.
  Translated in Child \cite{Ch20}.

\bibitem{Ge49} Gerhardt, C. (ed.)  \emph{Die philosophischen Schriften
  von Gottfried Wilhelm Leibniz}.  Berlin: Weidmann 1875--1890 (cited
  as [GP]).


\bibitem{Ge50} Gerhardt, C. (ed.)  \emph{Leibnizens mathematische
  Schriften}.  Berlin: Asher and Halle: H.W. Schmidt, 1849--63 (cited
  as [GM]).
%
%

\bibitem{Gr98} Grosholz, E.; Yakira, E.\, Leibniz's science of the
  rational.  Studia Leibnitiana Sonderheft, 26.  Franz Steiner Verlag
  Wiesbaden GmbH, Stuttgart, 1998.


\bibitem{Go08} Goldenbaum, U.; Jesseph, D. (Eds.)  \emph{Infinitesimal
  Differences.  Controversies between Leibniz and his Contemporaries}.
  Berlin New York: Walter de Gruyter, 2008.


\bibitem{He18} Hellman, G.; Shapiro, S.  \emph{Varieties of continua.
  From regions to points and back. Chapter 3 with {\O}ystein
  Linnebo}. Oxford: Oxford University Press, 2018.


\bibitem{Is90} Ishiguro, H.  \emph{Leibniz's philosophy of logic and
  language.}  Second edition.  Cambridge: Cambridge University Press,
  1990.





\bibitem{22b} Katz, M.; Kuhlemann, K.; Sherry, D.; Ugaglia, M.; van
  Atten, M.\, Two-track depictions of Leibniz's fictions.  \emph{The
    Mathematical Intelligencer} (2022).
  \url{https://doi.org/10.1007/s00283-021-10140-3},
  \url{https://arxiv.org/abs/2111.00922}


\bibitem{13f} Katz, M., Sherry, D.\, Leibniz's infinitesimals: Their
  fictionality, their modern implementations, and their foes from
  Berkeley to Russell and beyond.  \emph{Erkenntnis} \textbf{78}
  (2013), no.\;3, 571--625.  See
  \url{http://dx.doi.org/10.1007/s10670-012-9370-y} and
  \url{http://arxiv.org/abs/1205.0174}




\bibitem{Kn90} Knobloch, E.\, L'infini dans les math\'ematiques de
  Leibniz.  In \emph{L'infinito in Leibniz, Problemi e terminologia},
  Simposio Internazionale Roma, 6--8 novembre 1986, a cura di Antonio
  Lamarra, pp.\;33--51.  Roma: Edizioni dell' Ateneo, 1990.

\bibitem{Kn94} Knobloch, E.\, The infinite in Leibniz's mathematics -
  The historiographical method of comprehension in context.  In
  K.\;Gavroglu, J. Christianidis, E.\;Nicolaidis (eds.), Trends in the
  historiography of science, pp.\;265--278.  Dordrecht: Kluwer, 1994.


\bibitem{Kn99} Knobloch, E.\, Galileo and Leibniz: different
  approaches to infinity.  \emph{Archive for History of Exact
    Sciences} \textbf{54} (1999), no.\;2, 87--99.



\bibitem{La90} Lamarra, A.\, An anonymous criticism from Berlin to
  Leibniz's Philosophy: John Toland against Mathematical Abstractions.
  \emph{Studia Leibnitiana}, Sonderheft \textbf{16} (1990), 89--102.
%



\bibitem{Le69} Leibniz, G.W.\, \emph{Letter to Thomasius},
  20--30\;april\;1669. In Leibniz \cite{Le23} (A2.1$^2$. 23--38).


\bibitem{Le72} Leibniz, G.W.\, \emph{Confessio Philosophi}
  (1672/3). In Leibniz \cite{Le23} (A.VI. 3.  115--149).

\bibitem{Le76} Leibniz, G.W.\, \emph{Numeri Infiniti} (1676). In
  Leibniz \cite{Le23} (A.VI. 3.  496--504).


\bibitem{Le78} Leibniz, G.W.\, Letter to Elisabeth.  A.2.1$^2$.
  659--666, 1678.  Translated in \cite{Le89a}.


\bibitem{Le83} Leibniz. G.W.\, \emph{Elementa nova matheseos
  universalis} (1683). In Leibniz \cite{Le23} (A.VI. 4A. 513--524).


\bibitem{Le86} Leibniz, G.W.\, Meditatio nova de natura anguli
  contactus et osculi, horumque usu in practica mathesi ad figuras
  faciliores succedaneas difficilioribus substituendas.
  \emph{Act. Erudit. Lips.}  june 1686.  In Gerhardt \cite{Ge50}
  vol.\;VII, pp.\;326--329.  English translation at
  \url{http://www.17centurymaths.com/contents/Leibniz/ae4.pdf}



\bibitem{Le86b} Leibniz, G.W. \, \emph{Discours de Metaphysique}
  (1686).  In Leibniz \cite{Le23} (A.VI. 4B2. 1529--1589).

\bibitem {Le95a} Leibniz, G.W. \, \emph{Letter to l'Hospital},
  14/24\;june\;1695, in Gerhardt \cite{Ge50}, vol.\;I, pp.\;287--289.
%
%

\bibitem{Le95b} Leibniz, G.W. \, Responsio ad nonnullas difficultates
  a Dn.\;Bernardo Niewentiit circa methodum differentialem seu
  infinitesimalem motas.  \emph{Act. Erudit. Lips.} (1695).  In
  Gerhardt \cite{Ge50}, vol.\;V, pp.\;320--328.  A French translation
  by Parmentier is in \cite[pp.\;316--334]{Le1989}.
%
%


\bibitem{Le95c} Leibniz, G.W. \, \emph{Letter to Huygens},
  1\;july\;1695.  In Gerhardt \cite{Ge50}, vol.\;II, p.\; 205 and
  Briefwechsel p.\;757. See
  \url{https://dbnl.nl/tekst/huyg003oeuv10_01/huyg003oeuv10_01_0239.php}




\bibitem{Le95d} Leibniz, G. W.\, Specimen dynamicum, pro admirandis
  naturae legibus circa corporum vires et mutuas actiones detegendis,
  et ad suas causas revocandis.  \emph{Acta Eruditorum}, april 1695.
  In Gerhardt \cite{Ge50}, vol.\;VI, pp.\; 234--254.  Translated in
  Loemker \cite{Le89}, pp.\;435--444 (Part\;I).


\bibitem{Le01c} Leibniz, G.W.  \emph{Cum Prodiisset}\ldots{} mss Cum
  prodiisset atque increbuisset Analysis mea infinitesimalis\ldots{}
  (1701) in Gerhardt \cite{Ge46}, pp.\;39--50.  See
  \url{http://books.google.co.il/books?id=UOM3AAAAMAAJ}


\bibitem{Le02} Leibniz, G.W.\, \emph{Letter to Varignon},
  2\;february\;1702. In Gerhardt \cite{Ge50}, vol.\;IV, pp.\;91--95.

\bibitem{Le02b} Leibniz, G.W.\, \emph{Letter to Varignon},
  14\;april\;1702. In Gerhardt \cite{Ge50}, vol.\;IV, pp.\;97--99.

\bibitem{Le02c} Leibniz, G.W.\, \emph{Letter to Varignon},
  20\;june\;1702. In Gerhardt \cite{Ge50}, vol.\;IV, pp.\,106--110.


\bibitem{Le04} Leibniz, G.W.\, \emph{Letter to de Volder}, 30 june
  1704. In Gerhardt \cite{Ge49}, vol.\;II, pp.\;267--272.


\bibitem{Le06} Leibniz, G.W \, \emph{Letter to des Bosses},
  11\;march\;1706. In Gerhardt \cite{Ge49}, vol.\;II,
  pp.\;304--308. Translated in Look and Rutherford's Leibniz
  \cite{Le07}, pp.\;30--39.
%
%


\bibitem{Le16} Leibniz, G.W. \, \emph{Letter to Samuel Masson} (1716).
  In Gerhardt \cite{Ge49}, vol.\;VI, pp.\;624--629.  Translated in
  Ariew--Garber's Leibniz \cite{Le89a}, pp.\;225--230.

\bibitem{Le23} Leibniz, G.W.\, \emph{S\"amtliche Schriften und
  Briefe}, Darmstadt/Leipzig/Berlin: Akademie 1923-- (cited as
  A. series, volume, page number).  See
  \url{https://leibnizedition.de}

\bibitem{Le1989} Leibniz, G.W. \, \emph{La naissance du calcul
  diff\'erentiel.  26 articles des \emph{Acta Eruditorum}.}
  Translated from the Latin and with an introduction and notes by Marc
  Parmentier.  With a preface by Michel Serres.  Mathesis. Paris:
  Librairie Philosophique J. Vrin, 1989.  See
  \url{https://books.google.co.il/books?id=lfEy-OzaWkQC}


\bibitem{Le89} Leibniz, G.W. \emph{Philosophical papers and letters}.
  Second Edition.  Synthese Historical Library, Vol.\;2.\, Leroy
  E. Loemker, Editor and Translator.  Dordrecht--Boston--London:
  Kluwer Academic Publishers, 1989.


\bibitem{Le89a} Leibniz, G.W.\, \emph{Philosophical essays}.
  Translated and edited by Roger Ariew and Daniel Garber.
  Indianapolis: Hackett Publishing, 1989.  See
  \url{https://books.google.co.il/books?id=1xEeAt6FUl8C}


\bibitem{Le04b} Leibniz, G.W.\, \emph{Quadrature arithm\'etique du
  cercle, de l'ellipse et de l'hyperbole}.  Marc Parmentier
  (Trans. and Ed.)  /Latin text by Eberhard Knobloch (Ed.), Paris:
  J.\;Vrin, 2004.  See
  \url{https://books.google.co.il/books?id=fNTUULXHmQ0C}

\bibitem{Le04c} Leibniz, G.W.\, \emph{Confessio Philosophi: Papers
  Concerning the Problem of Evil, 1671--1678} (The Yale Leibniz
  Series), 2004, by G. W. Leibniz (Author), Robert C. Sleigh Jr.
  (Translator).  See
  \url{https://books.google.co.il/books?id=Jj3Z0dHp6_MC}


\bibitem{Le07} Leibniz, G.W.\, \emph{The Leibniz--Des Bosses
  Correspondence}.  Translated, edited, and with an Introduction by
  Brandon C. Look and Donald Rutherford.  New Haven: Yale University
  Press, 2007.


\bibitem{Le18} Leibniz, G.W.  \emph{Mathesis universalis}.  \'Ecrits
  sur la math\'ematique universelle.  Translated from the Latin and
  with an introduction and notes by David Rabouin.  Mathesis.  Paris:
  Librairie Philosophique J. Vrin, 2018.



\bibitem{Le98} Levey, S.\, Leibniz on mathematics and the actually
  infinite division of matter.  \emph{The Philosophical Review}
  \textbf{107} (1998), no.\,1, 49--96.


\bibitem{Le21} Levey, S.\, The continuum, the infinitely small, and
  the law of continuity in Leibniz.  In \emph{The History of Continua:
    Philosophical and Mathematical Perspectives}, pp.\,123--157.\,
  Edited by Stewart Shapiro and Geoffrey Hellman. Oxford: Oxford
  University Press, 2021.  See
  \url{https://doi.org/10.1093/OSO/9780198809647.003.0018}



\bibitem{Me94} 
%
%
Mercer, C.; Sleigh, R.\, Metaphysics: The early period to the
\emph{Discourse on Metaphysics}.\, In \emph{Cambridge Companion to
  Leibniz}, pp.\;67--123, Nicholas Jolley, Ed.\, Cambridge: Cambridge
University Press, 1994.




\bibitem{Kr} Murdoch, J. (1982). Infinity and continuity. In
  N. Kretzmann, A. Kenny, J.\;Pinborg, \& E.\;Stump (Eds.), \emph{The
    Cambridge History of Later Medieval Philosophy: From the
    Rediscovery of Aristotle to the Disintegration of Scholasticism,
    1100–1600},\, pp.\;564--592.  Cambridge: Cambridge University
  Press.


\bibitem{Ra13} Rabouin, D.\, ``Analytica Generalissima Humanorum
  Cognitionum''.  Some reflections on the relationship between logical
  and mathematical analysis in Leibniz.  \emph{Studia Leibnitiana},
  2013, Bd.\;\textbf{45}, H.\,1 (2013), pp.\,109--130.


\bibitem{Ra15} Rabouin, D.\, Leibniz's rigorous foundations of the
  method of indivisibles Or How to Reason with Impossible Notions.  In
  Vincent Jullien (Ed.), \emph{Seventeenth-Century Indivisibles
    Revisited}, Science Networks.  Historical Studies,
  Vol.\;49. Basel: Birkh\"auser, 2015, pp.\;347--364.


\bibitem{Ra15b} Rabouin, D.\, The difficulty of being simple: On some
  interactions between mathematics and philosophy in Leibniz's
  analysis of notions.  In N. B. Goethe et al. (eds.),
  \emph{G. W. Leibniz, Interrelations between Mathematics and
    Philosophy}, Archimedes 41, Netherlands: Springer, 2015.


\bibitem{Ra21} Rabouin, D.\, Can one be a fictionalist and a platonist
  at the same time?  Lessons from Leibniz (2021).



\bibitem{Ra20} Rabouin, D.; Arthur, R.\, Leibniz's syncategorematic
  infinitesimals II: their existence, their use and their role in the
  justification of the differential calculus.  \emph{Archive for
    History of Exact Sciences} \textbf{74} (2020), 401--443.



\bibitem{Sc05} Schubring, G.  \emph{Conflicts between generalization,
  rigor, and intuition.  Number concepts underlying the development of
  analysis in 17--19th Century France and Germany}.\, Sources and
  Studies in the History of Mathematics and Physical Sciences. New
  York: Springer-Verlag, 2005.


\bibitem{14c} Sherry, D.; Katz, M.\, Infinitesimals, imaginaries,
  ideals, and fictions.  \emph{Studia Leibnitiana} \textbf{44} (2012),
  no.\;2, 166--192.  See \url{http://www.jstor.org/stable/43695539}
  and \url{https://arxiv.org/abs/1304.2137} (Article was published in
  2014 even though the journal issue lists the year as 2012)


\bibitem{Si99} Sierksma, G.; Sierksma, W.\, The great leap to the
  infinitely small.  Johann Bernoulli: mathematician and philosopher.
  \emph{Ann. of Sci.}  \textbf{56} (1999), no.\;4, 433--449.


\bibitem{Th82} Thomason, S.\, Euclidean infinitesimals.  \emph{Pacific
  Philosophical Quarterly} \textbf{63} (1982), 168--185.


\bibitem{To11} Toland, J.\, Remarques Critiques sur le Syst\^eme de
  Mr.\ Leibnitz de l'Harmonie pr\'e\'etablie; o\`u l'on recherche, en
  passant, pourquoi les Syst\^emes Metaphysiques des Mathematiciens
  ont moins de clart\'e que ceux des autres: \'ecrites par ordre de sa
  Majest\'e la feu\"e Reine de Prusse.  \emph{Histoire Critique de la
    Republique des Lettres}, Tome \textbf{XI} (1716), 115--133
  (published anonymously).


\bibitem{Ug22} Ugaglia, M.\, Possibility vs Iterativity: Leibniz and
  Aristotle on the Infinite. In \emph{Thinking and Calculating. Essays
    on Logic, its History and its Applications in honor of Massimo
    Mugnai}, Francesco Ademollo, Fabrizio Amerini, Vincenzo De Risi
  (Eds.), Berlin, Springer 2022.



\bibitem{Wo98} Woolhouse, R.\, John Toland and `Remarques Critiques
  sur le Syst\^eme de Monsr.\ Leibnitz de l'Harmonie pr\'e\'etablie'.
  \emph{The Leibniz Review} \textbf{8} (1998), 80--87.  See
  \url{https://doi.org/10.5840/leibniz199883}


\bibitem{Wo97} Woolhouse, R.; Francks, R. (Eds.)  \emph{Leibniz's `New
  System' and Associated Contemporary Texts}.  Oxford: Clarendon
  Press, 1997.


\end{thebibliography}
\end{document}